\documentclass[a4paper,12pt,twoside]{article}
\usepackage[francais]{babel}
\usepackage{amssymb, amsmath, eucal}
\usepackage[all]{xy}
\usepackage{mathrsfs}
\usepackage{leqno}
\begin{document}
\textwidth15.5cm
\textheight22.5cm
\voffset=-13mm
\newtheorem{The}{Theorem}[subsection]
\newtheorem{Lem}[The]{Lemma}
\newtheorem{Prop}[The]{Proposition}
\newtheorem{Cor}[The]{Corollary}
\newtheorem{Rem}[The]{Remark}
\newtheorem{Titre}[The]{\!\!\!\! }
\newtheorem{Conj}[The]{Conjecture}
\newtheorem{Question}[The]{Question}
\newtheorem{Prob}[The]{Problem}
\newtheorem{Def}[The]{Definition}
\newtheorem{Not}[The]{Notation}
\newtheorem{Ex}[The]{Example}
\newcommand{\C}{\mathbb{C}}
\newcommand{\R}{\mathbb{R}}
\newcommand{\N}{\mathbb{N}}
\newcommand{\Q}{\mathbb{Q}}

\begin{center}

{\Large\bf Effective Local Finite Generation of Multiplier Ideal Sheaves\footnote{2000 Mathematics Subject Classification: 32C35, 32U05, 32A36}}

\end{center}

\begin{center}

{\large\bf Dan Popovici}

\end{center}

\vspace{1ex}

\noindent {\small {\bf Abstract.} Let $\varphi$ be a psh function on a bounded pseudoconvex open set $\Omega \subset \C^n$, and let ${\cal I}(\varphi)$ be the associated multiplier ideal sheaf. Motivated by resolution of singularities issues, we establish an effective version of the coherence property of ${\cal I}(m\varphi)$ as $m\rightarrow +\infty$. Namely, we estimate the order of growth in $m$ of the number of generators needed to engender ${\cal I}(m\varphi)$ on a fixed compact subset, as well as the growth of the coefficients featuring in the decomposition of local sections as linear combinations over ${\cal O}_{\Omega}$ of finitely many generators. The main idea is to use Toeplitz concentration operators involving Bergman kernels associated with singular weights. Our approach relies on asymptotic integral estimates of singularly weighted Bergman kernels of independent interest. In the second part of the paper, we estimate the additivity defect of multiplier ideal sheaves already known to be subadditive by a result of Demailly, Ein, and Lazarsfeld. This implies that the decay rate of ${\cal I}(m\varphi)$ is not far from being linear if the singularities of $\varphi$ are reasonable.}

\vspace{3ex}

\subsection{Introduction}\label{subsection:intro}

Let $X$ be a complex manifold of complex dimension $n$, and let $\varphi$ be a plurisubharmonic (psh) function on $X$. Following Nadel ([Nad90]), one can associate with $\varphi$ the multiplier ideal sheaf ${\cal I}(\varphi)\subset {\cal O}_X$ whose fibre at a point $x\in \Omega$ is defined as the set of germs of holomorphic functions $f\in {\cal O}_{\Omega, \, x}$ such that $|f|^2 \, e^{-2\varphi}$ is integrable with respect to the Lebesgue measure in some local coordinates in a neighbourhood of $x$. The main purpose of the present paper is to study the variation of ${\cal I}(m\varphi)$ as $m\rightarrow +\infty$.

 The sequence of multiplier ideal sheaves ${\cal I}(m\varphi)$ is easily seen to be nonincreasing as $m\rightarrow +\infty$. Indeed, to offset the possible nonintegrability of $e^{-2m\varphi}$ near points $x\in \Omega$ where $\varphi(x)=-\infty$, holomorphic germs $f\in {\cal I}(m\varphi)_x$ may need to vanish to increasingly high orders at $x$ as $m$ increases. On the other hand, the subadditivity property of multiplier ideal sheaves established by Demailly, Ein, and Lazarsfeld (cf. [DEL00]) implies that \\

\hspace{20ex} ${\cal I}(m\varphi)\subset {\cal I}(\varphi)^m$, \\

\noindent while the inclusion is strict in general. In other words, ${\cal I}(m\varphi)$ may decrease more quickly than linearly as $m\rightarrow +\infty$. The main thrust of the ensuing development is to obtain an effective control on the decay of ${\cal I}(m\varphi)$ as $m\rightarrow +\infty$. Since the problems dealt with throughout this paper are local in nature, we shall focus our attention on psh functions $\varphi$ defined on a bounded pseudoconvex open set $\Omega \subset \mathbb{\C}^n$. Examples of such functions are provided by the so-called psh functions with {\bf analytic singularities}, namely those functions that can be locally written as \\

\hspace{5ex}  $(\star)$  \hspace{15ex} $\displaystyle \varphi = \frac{c}{2} \, \log (|g_1|^2 + \dots + |g_N|^2)$, \\

\noindent for some holomorphic functions $g_1, \dots , \, g_N$ on $\Omega$, and some constant $c>0$. The $-\infty$-poles (or singularities) of $\varphi$ are precisely the common zeroes of $g_1, \dots , \, g_N$, and they only depend, up to an equivalence of singularities, on the sheaf generated by $g_1, \dots , \, g_N$.  \\

 The issue of the decaying rate of ${\cal I}(m\varphi)$ is addressed in two ways. \\

  First, let ${\cal H}_{\Omega}(m\varphi)$ be the Hilbert space of holomorphic functions $f$ on $\Omega$ such that $|f|^2 \, e^{-2m\varphi}$ is integrable with respect to the Lebesgue measure on $\Omega$. It is well-known that the ideal sheaf ${\cal I}(m\varphi)$ is coherent and generated as an ${\cal O}_{\Omega}$-module by an arbitrary orthonormal basis $(\sigma_{m, \, j})_{j\in \mathbb{N}^{\star}}$ of ${\cal H}_{\Omega}(m\varphi)$ (see e.g. [Dem93, Lemma 4.4]). By the strong Noetherian property of coherent sheaves, it is then generated, on every relatively compact open subset $\Omega' \subset\subset \Omega$, by only finitely many $\sigma_{m, \, j}$'s. Our first goal is to make this local finite generation property effective; in other words, to estimate the number $N_m$ of generators needed, as well as the growth rate of the (holomorphic function) coefficients appearing in the decomposition of an arbitrary section of ${\cal I}(m\varphi)$ on $\Omega'$ as a finite linear combination of $\sigma_{m, \, j}$'s, as $m\rightarrow +\infty$. The first set of results can be summed up as follows.

 \begin{The}\label{The:introd1}Let $\varphi$ be a strictly psh function on $\Omega\subset \C^n$ such that $i\partial\bar{\partial}\varphi \geq C_0\, \omega$ for some constant $C_0>0$. Let $B:=B(x, \, r) \subset\subset\Omega$ be an arbitrary open ball. Then, there exist a ball $B(x, \, r_0)\subset\subset B(x, \, r)$ and $m_0=m_0(C_0)\in \mathbb{N}$, such that for every $m\geq m_0$ the following property holds. Every $g\in {\cal H}_B(m\varphi)$ admits, with respect to some suitable finitely many elements $\sigma_{m, \, 1}, \dots \sigma_{m, \, N_m}$ in a suitable orthonormal basis $(\sigma_{m, \, j})_{j\in \N^{\star}}$ of ${\cal H}_{\Omega}(m\varphi)$, a decomposition:  \\

$\displaystyle g(z) = \sum\limits_{j=1}^{N_m} b_{m, \, j}(z)\, \sigma_{m, \, j}(z)$,  \hspace{3ex}  $z\in B(x, \, r_0),$  \\

\noindent with some holomorphic functions $b_{m, \, j}$ on $B(x, \, r_0),$ satisfying: \\

 $\displaystyle\sup\limits_{B(x, \, r_0)}\sum\limits_{j=1}^{N_m}|b_{m, \, j}|^2 \leq C \, N_m \, \int_B|g|^2\, e^{-2m\varphi}<+\infty,$ \\

\noindent where $C>0$ is a constant depending only on $n, \, r,$ and the diameter of $\Omega.$ \\

\noindent Moreover, if $\varphi$ has analytic singularities, then $N_m\leq C_{\varphi} \, m^n$ for $m>> 1$, where $C_{\varphi}>0$ is a constant depending only on $\varphi$, $B$, and $n$. 

\end{The}

This can be seen as a local counterpart to the effective version of the global generation of multiplier ideal sheaves proved in [Siu02, Theorem 2.1] as a step towards the invariance of plurigenera for projective manifolds not necessarily of general type. Unlike in Siu's case, where spaces of global sections over compact manifolds were already finite dimensional, the corresponding spaces ${\cal H}_{\Omega}(m\varphi)$ of local sections are infinite dimensional. To achieve a finite number of generators, a restriction to a compact subset is necessary. This involves appropriate choices of generators with estimates obtained from a study of Bergman kernels. The main technique is a concentration procedure of $L^2$ norms on compact subsets carried out by means of Toeplitz operators.

Second, to complement the subadditivity property of Demailly, Ein, and Lazarsfeld, we obtain a superadditivity result showing that the variation of ${\cal I}(m\varphi)$, though not necessarily linear in $m$, is not far from being so if the singularities of $\varphi$ are analytic. Specifically, the following holds.

  \begin{The}\label{The:introd2} Let $\varphi$ be psh with analytic singularities on $\Omega \subset \C^n.$ Then, for every $\Omega'\subset\subset \Omega$ and every $\delta >0,$ we have:  \\

$\noindent {\cal I}(m\varphi)^p_{|\Omega'} \subset {\cal I}(mp(1-\delta)\varphi)_{|\Omega'},$ \hspace{3ex}  for every $m\geq \frac{n+2}{c\, \delta},$ and $p\in \mathbb{N}^{\star}.$

\end{The}

 Here is an outline of our approach. In section \ref{subsection:Toeplitz}, we recall the definition of Toeplitz concentration operators and observe some basic properties. We then extend to the case of analytic singularities some asymptotic integral estimates on the Bergman kernel previously obtained for smooth weight functions by N. Lindholm ([Lin01]) and B. Berndtsson ([Ber03]). While these results may have an interest of their own, we apply them in our case at hand to estimate the number $N_m$ of local generators needed for ${\cal I}(m\varphi)$. In section \ref{subsection:generation}, H\"ormander's $L^2$ estimates ([H\"or65]) and Skoda's $L^2$ division theorem ([Sko72]) are used to prove an effective local finite generation property of ${\cal I}(m\varphi)$, and thus complete the proof of Theorem \ref{The:introd1}. Section \ref{subsection:regularization} deals with the different but related issue of estimating the subadditivity defect of multiplier ideal sheaves. Demailly's procedure for locally approximating arbitrary psh functions by psh functions with analytic singularities ([Dem92], Proposition 3.1) is revisited to get improved effective estimates on fixed-sized subsets. The key ingredient is once more Skoda's $L^2$ division theorem. Theorem \ref{The:introd2} then follows as a corollary.

 We believe that such effective estimates on multiplier ideal sheaves could be useful in the birational geometry of complex manifolds when singularities are resolved by blowing up the multiplier ideal sheaves encoding them.

\vspace{3ex}

\subsection {Toeplitz concentration operators}\label{subsection:Toeplitz} 

 This section is intended to lay the groundwork for the proof of Theorem \ref{The:introd1} by hinting at how a suitable orthonormal basis of ${\cal H}_{\Omega}(m\varphi)$ and finitely many local generators of ${\cal I}(m\varphi)$ will be chosen. The actual choices will be spelt out in the next section. The key point is an estimate of the growth order of the number $N_m$ of local generators via asymptotic estimates on the Bergman kernel of independent interest.

Let $B\subset\subset \Omega$ be a relatively compact open subset, typically a ball. Following [Lin01] and [Ber03] (themselves inspired by [Lan67]), we consider, for every $m\in \mathbb{N}^{\star}$, the following Toeplitz concentration operator on $B$: \\

$T_{B, \, m} : {\cal H}_{\Omega}(m\varphi) \rightarrow {\cal H}_{\Omega}(m\varphi)$, \hspace{3ex}  $T_{B, \, m}(f)=P_m (\chi_B \, f)$,  \\

\noindent where $\chi_B$ is the characteristic function of $B$, and $P_m:L^2(\Omega, \, e^{-2m\varphi}) \rightarrow {\cal H}_{\Omega}(m\varphi)$ is the orthogonal projection from the Hilbert space of (equivalence classes of) measurable functions $f$ for which $|f|^2 \, e^{-2m\varphi}$ is Lebesgue integrable on $\Omega$, onto the closed subspace of holomorphic such functions. It is easy to see that \\

$T_{B, \, m}(f)= \chi_B \, f -u$, \\

\noindent where $u$ is the solution of the equation $\bar{\partial}u = \bar{\partial}(\chi_B \, f)$ of minimal $e^{-2m\varphi}$-weighted $L^2$-norm. Alternatively, if we consider the Bergman kernel: \\

$K_{m\varphi} : \Omega\times \Omega \rightarrow \mathbb{\C}$, \hspace{3ex}  $\displaystyle K_{m\varphi}(z, \, \zeta) = \sum\limits_{j=1}^{+\infty} \sigma_{m, \, j}(z) \, \overline{\sigma_{m, \, j}(\zeta)}$,  \\

\noindent its reproducing property shows that the concentration operator is also given as \\

$\displaystyle T_{B, \, m}(f)(z) = \int_B K_{m\varphi}(z, \, \zeta) \, f(\zeta) \, e^{-2m\varphi(\zeta)}\, d\, \lambda(\zeta)$,  \hspace{3ex}  $z\in \Omega$,  \\

\noindent where $d\, \lambda$ is the Lebesgue measure. It is clear that $T_{B, \, m}$ is a compact operator for it is defined by a square integrable kernel. Its eigenvalues $\lambda_{m, \, 1} \geq \lambda_{m, \, 2} \geq \dots $ lie in the open interval $(0, \, 1)$. If $f\in {\cal H}_{\Omega}(m\varphi)$ is an eigenvector of $T_{B, \, m}$ corresponding to some eigenvalue $\lambda$, we see that $\langle\langle T_{B, \, m}(f), \, f \rangle\rangle = \lambda \, ||f||^2$, and implicitly:  \\

$\displaystyle \lambda = \displaystyle \frac{\int_B |f|^2 \, e^{-2m\varphi}}{\int_{\Omega} |f|^2 \, e^{-2m\varphi}}$. \\

\noindent Therefore, $f\in {\cal H}_{\Omega}(m\varphi)$ is an eigenvector of $T_{B, \, m}$ if and only if \\

$\displaystyle f(z) = \frac{\int_{\Omega} |f|^2 \, e^{-2m\varphi}}{\int_B |f|^2 \, e^{-2m\varphi}}\, \int_B K_{m\varphi}(z, \, \zeta) \, f(\zeta) \, e^{-2m\varphi(\zeta)}\, d\, \lambda(\zeta)$, \hspace{3ex} $z\in \Omega$,  \\

\noindent which amounts to having:  \\

$f(z) = \displaystyle \frac{\int_{\Omega} |f|^2 \, e^{-2m\varphi}}{\int_B |f|^2 \, e^{-2m\varphi}}\, \sum\limits_{j=1}^{+\infty} \sigma_{m, \, j}(z) \, \int_B f(\zeta)\, \overline{\sigma_{m, \, j}(\zeta)} \, e^{-2m\varphi(\zeta)}\, d\, \lambda(\zeta)$, \hspace{1ex} $z\in \Omega$.  \\

\noindent On the other hand, every $f\in {\cal H}_{\Omega}(m\varphi)$ has a Hilbert space decomposition with Fourier coefficients as:  \\

$f(z) = \displaystyle \sum\limits_{j=1}^{+\infty} \bigg (\int_{\Omega}f(\zeta)\, \overline{\sigma_{m, \, j}(\zeta)} \, e^{-2m\varphi(\zeta)}\, d\, \lambda(\zeta) \bigg )\, \sigma_ {m, \,j}(z)$,  \hspace{3ex} $z\in \Omega$.  \\

\noindent The uniqueness of the decomposition into a linear combination of elements in an orthonormal basis implies the following simple observation.

\begin{Lem}\label{Lem:concentratedbasis} A function $f\in {\cal H}_{\Omega}(m\varphi)$ is an eigenvector of $T_{B, \, m}$ if and only if  \\

$\displaystyle \frac{1}{\int_{\Omega}|f|^2 \, e^{-2m\varphi}} \, \int_{\Omega} f\, \bar{\sigma}_{m, \, j} \, e^{-2m\varphi} = \frac{1}{\int_{B}|f|^2 \, e^{-2m\varphi}} \, \int_{B} f\, \bar{\sigma}_{m, \, j} \, e^{-2m\varphi}$,  

\vspace{3ex}

\noindent for every $j\geq 1$.

\end{Lem}

\noindent We will be now studying the behaviour of the eigenvalues of $T_{B, \, m}$ as $m\rightarrow +\infty$.   Let us fix an orthonormal basis $(\sigma_{m, \, j})_{j\in \mathbb{N}^{\star}}$ of ${\cal H}_{\Omega}(m\varphi)$ made up of eigenvectors of $T_{B, \, m}$ corresponding respectively to its eigenvalues $\lambda_{m, \, 1} \geq \lambda_{m, \, 2} \geq \dots $ listed nonincreasingly. Let $\varepsilon > 0$. Since $T_{B, \, m}$ is a compact operator, there are (if any) at most finitely many eigenvalues $\lambda_{m, \, 1} \geq \lambda_{m, \, 2} \geq \dots \geq \lambda_{m, \, N_m} \geq 1-\varepsilon$. In other words, \\

$\displaystyle \int_B |\sigma_{m, \, 1}|^2 \, e^{-2m\varphi} \geq \dots \geq \int_B |\sigma_{m, \, N_m}|^2 \, e^{-2m\varphi} \geq 1-\varepsilon > \int_B |\sigma_{m, \, k}|^2 \, e^{-2m\varphi}$, \\

\noindent for every $k\geq N_m +1$.  \\

Lemma \ref{Lem:concentratedbasis} above shows, in particular, that the restrictions to $B$ of the $\sigma_{m, \, j}$'s are still orthogonal to one another. If we let 

\begin{equation}\label{eqn:defBergman}\displaystyle B_{m\varphi}(z):=K_{m\varphi}(z, \, z)=\sum\limits_{j=1}^{+\infty}|\sigma_{m, \, j}(z)|^2 = \sup\limits_{f\in B(1)}|f(z)|^2,  \hspace{3ex} z\in \Omega,\end{equation}

\noindent where $B(1)$ is the unit ball in ${\cal H}_{\Omega}(m\varphi)$, the traces of $T_{B, \, m}$ and $T_{B, \, m}^2$ are easily computed as: \\

 $\displaystyle\mbox{Tr}\, (T_{B, \, m})=\int_B B_{m\varphi}(z)\, e^{-2m\varphi(z)} \, d\, \lambda(z)$, \hspace{3ex} and \\

 $\displaystyle\mbox{Tr}\, (T_{B, \, m}^2)= \int_{B\times B} |K_{m\varphi}(z, \, \zeta)|^2\, e^{-2m\varphi(z)}\, e^{-2m\varphi(\zeta)} \, d\, \lambda(z)\, d\, \lambda(\zeta)$. \\

\noindent We are thus naturally led to a search for asymptotic estimates of the Bergman kernel, preferably in the general context of possibly singular psh functions $\varphi$. Bearing in mind the analogy with the $L^2$ volume of a line bundle defined on a compact complex manifold in terms of global sections and characterized in terms of curvature currents, we set the following.

\begin{Def}\label{Def:volume} The {\it volume} on $B\subset\subset \Omega$ of a psh function $\varphi$ on $\Omega$ is defined as \\

$\displaystyle v_B(\varphi):= \limsup\limits_{m\rightarrow +\infty} \frac{n !}{m^n}\, \int_B B_{m\varphi}\, e^{-2m\varphi}\, d\, \lambda$, \\

\noindent where $d\, \lambda$ is the Lebesgue measure.

\end{Def}

It is our intention to study for which psh functions $\varphi$ the volume $v_B(\varphi)$ is finite. The outcome will be a control of the number $N_m$ of eigenvalues of $T_{B, \, m}$ which are $\geq 1-\varepsilon$. These eigenvalues will be seen to correspond to the local generators of ${\cal I}(m\varphi)$ whose growth rate is estimated in Theorem \ref{The:introd1}. The motivation lies in the following asymptotic estimates for the Bergman kernel associated with a smooth $\varphi$, due to N. Lindholm (cf. Theorems 10, 11, and 13 in [Lin01]), and subsequently rewritten in a slightly different setting by B.Berndtsson (cf. Theorems 2.3, 2.4, and 3.1 in [Ber03]). The standard K\"ahler form on $\C^n$ will be denoted throughout by $\omega$.

\begin{The}\label{The:Bo} (Lindholm, Berndtsson) Let $\varphi$ be a $C^{\infty}$ psh function such that $i\partial\bar{\partial}\varphi \geq C_0 \, \omega$ on $\Omega\subset \subset \C^n$ for some constant $C_0>0$. Then:  \\

\vspace{1ex}

\noindent $(a)$ \,  $\displaystyle v_B(\varphi) = \lim\limits_{m\rightarrow +\infty} \frac{n!}{m^n} \, \int_B B_{m\varphi}(z)\, e^{-2m\varphi(z)} \, d\, \lambda(z) = \frac{2^n}{\pi^n} \, \int_B (i\partial\bar{\partial}\varphi)^n$;  \\

\vspace{1ex} 

\noindent $(b)$ \, the sequence of measures on $\Omega \times \Omega$ defined as  \\

\hspace{3ex} $\displaystyle\frac{n!}{m^n} \, |K_{m\varphi}(z, \, \zeta)|^2\, e^{-2m\varphi(z)}\, e^{-2m\varphi(\zeta)} \, d\, \lambda(z) \, d\, \lambda(\zeta)$ \\

\noindent  converges to $\displaystyle \frac{2^n}{\pi^n} \, (i\partial\bar{\partial}\varphi)^n \wedge [\Delta ]$ in the weak topology of measures, as $m\rightarrow +\infty$, where $[\Delta ]$ is the current of integration on the diagonal of $\Omega \times \Omega$; \\

\noindent $(c)$ \,  finally, $\displaystyle \lim\limits_{m\rightarrow +\infty} \frac{n!}{m^n} \, \sum\limits_{j=1}^{+\infty}\lambda_{m, \, j} = \lim\limits_{m\rightarrow +\infty} \frac{n!}{m^n} \, \sum\limits_{j=1}^{+\infty}\lambda_{m, \, j}^2 =\frac{2^n}{\pi^n} \, \int_B (i\partial\bar{\partial}\varphi)^n$,    \\

\noindent and since $\lambda_{m, \, 1} \geq \dots \geq \lambda_{m, \, N_m} \geq 1-\varepsilon  > \lambda_{m, \, k}$, for $k\geq N_m+1$, we get:  \\

$\displaystyle \lim\limits_{m\rightarrow +\infty} \frac{n!}{m^n} \, N_m = \frac{2^n}{\pi^n} \, \int_B (i\partial\bar{\partial}\varphi)^n$ \hspace{2ex} and \hspace{2ex} $\displaystyle \lim\limits_{m\rightarrow +\infty} \frac{n!}{m^n}\sum\limits_{k=N_m + 1}^{+\infty}\lambda_{m, \, k} = 0$.

\end{The}

 We now single out the major steps in the proof given to this theorem in [Ber03, p. 4-6], as they will be subsequently adapted to a more general context. The proof of $(a)$ hinges on two facts. First, if $\lambda_1, \dots , \, \lambda_n$ are the eigenvalues of $i\partial\bar{\partial}\varphi$ with respect to $\omega$, the following pointwise convergence is established:

\begin{equation}\label{eqn:pointwiseestimate}\displaystyle \lim\limits_{m\rightarrow +\infty} \frac{1}{m^n}\, B_{m\varphi}(z)\, e^{-2m\varphi(z)} = \frac{2^n}{\pi^n} \, \lambda_1(z) \dots \lambda_n(z), \hspace{3ex} z\in \Omega.\end{equation} 

\noindent Indeed, having fixed an arbitrary point $x\in \Omega$, the quadratic part $\varphi_0$ in the Taylor expansion of $\varphi$ near $x$ is diagonalized in some local coordinates centred at $x$. While the distance between $\varphi$ and $\varphi_0$ is under control, the mean-value inequality is applied to $\varphi_0$ to get the upper estimate of the left-hand side by the right-hand side. The converse estimate is obtained by means of H\"ormander's $L^2$ estimates allowing one to construct elements in ${\cal H}_{\Omega}(m\varphi)$ with a prescribed value at $x$ and a global $L^2$ norm under control. These are local procedures carried out on small balls centred at $x$, and the desired estimates are obtained in the limit while shrinking the balls to $x$.

 Second, the mean-value inequality argument alluded to above also gives the following uniform estimate on a relatively compact open subset:

\begin{equation}\label{eqn:unifestimate}  0 < \frac{n!}{m^n}\, B_{m\varphi}(z)\, e^{-2m\varphi(z)} \leq C \, \frac{2^n}{\pi^n} \, \lambda_1(z) \dots \lambda_n(z), \hspace{3ex} z\in B,\end{equation}

\noindent if $m$ is big enough, for a constant $C>0$ independent of $m$. This is possible since $\varphi$ is $C^{\infty}$ and $B$ is relatively compact. One can then conclude by dominated convergence.

 We will now prove that Theorem \ref{The:Bo} is still essentially valid if we allow analytic singularities for $\varphi$ (see definition $(\star)$ in the introduction), provided the current $i\partial\bar{\partial}\varphi$ is replaced throughout by its absoulutely continuous part $(i\partial\bar{\partial}\varphi)_{ac}$ with respect to the Lebesgue measure in the Lebesgue decomposition of its measure coefficients (into an absolutely continuous part and a singular part). We clearly have $(i\partial\bar{\partial}\varphi)_{ac}=i\partial\bar{\partial}\varphi$ if $\varphi$ is $C^{\infty}$. 

 We will split the analysis of the analytic singularity case into two steps according to whether the coefficient of these singularities is an integer or not. \\

\noindent {\bf $(a)$ \,\, Analytic singularities with an integer coefficient}\\

 In this case, we have a complete analogue of Theorem \ref{The:Bo}. 

\begin{The}\label{The:integercoeff} Let $\varphi = \frac{p}{2} \, \log(|g_1|^2 + \dots + |g_N|^2) + u$ for some holomorphic functions  $g_1, \dots , \, g_N$, some $p\in \N$, and some $C^{\infty}$ function $u$ on $\Omega\subset\subset \C^n$. Assume $i\partial\bar{\partial}\varphi \geq C_0 \, \omega$ for some constant $C_0>0$. Then, for $B\subset\subset \Omega$, we have: \\

\noindent $(a)$\, $\displaystyle v_B(\varphi)=\lim\limits_{m\rightarrow +\infty} \frac{n!}{m^n}\, \int\limits_{B}B_{m\varphi}\, e^{-2m\varphi}\, d\, \lambda = \frac{2^n}{\pi^n} \, \int\limits_{B} (i\partial\bar{\partial}\varphi)_{ac}^n<+\infty.$ \\

\noindent $(b)$\,  the analogue of $(b)$ in Theorem \ref{The:Bo} holds, with convergence to \\

$\displaystyle \frac{2^n}{\pi^n}\, (i\partial\bar{\partial}\varphi)^n_{ac}\wedge [\Delta]$ \hspace{3ex} in the weak topology of measures; \\

\noindent $(c)$ \, the analogue of $(c)$ in Theorem \ref{The:Bo} holds, with \\

$\displaystyle \lim\limits_{m\rightarrow +\infty} \frac{n!}{m^n} \, N_m = \frac{2^n}{\pi^n} \, \int_B (i\partial\bar{\partial}\varphi)^n_{ac} < +\infty$.

\end{The}

 The key observation in this new setting is that $(i\partial\bar{\partial}\varphi)^n_{ac}$ is of locally finite mass, and thus the integrals involving $(i\partial\bar{\partial}\varphi)_{ac}^n$ above are finite. Indeed, we can resolve the analytic singularities of $\varphi$ by blowing up the ideal sheaf ${\cal I}$ generated as an ${\cal O}_{\Omega}$-module by $g_1, \dots , \, g_N$. According to Hironaka, there exists a proper modification $\mu:\tilde{\Omega}\rightarrow \Omega$, arising as a locally finite sequence of smooth-centred blow-ups, such that $\mu^{\star} {\cal I}={\cal O}(-D)$ for an effective divisor $D$ on $\tilde{\Omega}$. We thus get the following Zariski decomposition of the pull-back current: \\

\hspace{6ex} $\displaystyle\mu^{\star}(i\partial\bar{\partial}\varphi)=\alpha + c\, [D]$, \hspace{6ex} on $\tilde{\Omega}$, \\

\noindent with a $C^{\infty}$ closed $(1, \, 1)$-form $\alpha \geq 0$, where $[D]$ stands for the $D$-supported current of integration on $D$. If \\

\hspace{6ex} $V:=\{g_1 = \dots = g_N=0\}$ \\

\noindent is the singular set of $\varphi$, we clearly have $\int_B (i\partial\bar{\partial}\varphi)^n_{ac}=\int_{\mu^{-1}(B)}\alpha^n$, and this quantity is finite since the smooth volume form $\alpha^n$ has locally finite mass. If $\lambda_1, \dots , \, \lambda_n$ are the eigenvalues of $i\partial\bar{\partial}\varphi$ with respect to $\omega$ on $\Omega\setminus V$, this means that the product $\lambda_1 \dots \lambda_n$ is integrable on $B\setminus V$. \\

\noindent {\it Proof of Theorem \ref{The:integercoeff}.} The overall idea is to run Berndtsson's proof of the smooth case (cf. [Ber03], p. 4-6) on $\Omega \setminus V$ where $\varphi$ is $C^{\infty}$. Thanks to the local nature of its proof, the pointwise convergence (\ref{eqn:pointwiseestimate}) still holds at points $x\in \Omega \setminus V$. The main difficulty stems from the uniform estimate (\ref{eqn:unifestimate}) not being immediately clear near the singular set $V$. We claim, however, that the analogous uniform upper estimate does hold outside the singular set, namely:

\begin{equation}\label{eqn:unifsingestimate}  0 < \frac{n!}{m^n}\, B_{m\varphi}(z)\, e^{-2m\varphi(z)} \leq C \, \frac{2^n}{\pi^n} \, \lambda_1(z) \dots \lambda_n(z), \hspace{3ex} z\in B \setminus V,\end{equation}

\noindent if $m$ is big enough, for a constant $C>0$ independent of $m$.

Assuming this for the moment, the strong Noetherian property satisfied by the coherent sheaf ${\cal I}(m\varphi)$ implies that $B_m\, e^{-2m\varphi}$ is integrable on $B\subset\subset \Omega$, since it is dominated there by a constant multiplied by a finite sum $\sum|\sigma_{m, \, j}|^2\, e^{-2m\varphi}$. Its integral on $B$ equals therefore its integral on $B\setminus V$ (as $V$ is Lebesgue-negligible). Using now the key observation that $\lambda_1(z) \dots \lambda_n(z)$ is integrable on $B\setminus V$, $(a)$ follows by dominated convergence as in the smooth case.

The proof of $(b)$ in [Ber03] (p. 6-7) can be repeated on $\Omega \setminus V$ and extended across $V$ in a similar way. Explicitly, what we have to prove is that for every compactly supported continuous function $g$ on $\Omega \times \Omega$, we have:  \\

$\displaystyle\lim\limits_{m\rightarrow +\infty} \frac{1}{m^n} \int\limits_{\Omega \times \Omega} g(z, \, \zeta)\, |K_{m\varphi}(z, \, \zeta)|^2\, e^{-2m\varphi(z)}\, e^{-2m\varphi(\zeta)} \, d\, \lambda(z) \, d\, \lambda(\zeta) = $  \\

\hspace{40ex} $\displaystyle = \frac{2^n}{\pi^n\, n!} \int\limits_{\Omega \setminus V} g(z, \, z)\, (i\partial\bar{\partial}\varphi)^n(z).$  \\

\noindent Again, $(i\partial\bar{\partial}\varphi)_{ac}^n$ having locally finite mass on $\Omega$ implies the well-definedness of $(i\partial\bar{\partial}\varphi)_{ac}^n \wedge [\Delta]$ as a complex measure on $\Omega \times \Omega$. If $M:=\sup|g|$, we notice that:  \\

$\bigg|\displaystyle  \frac{1}{m^n} \int\limits_{\Omega} g(z, \, \zeta)\, |K_{m\varphi}(z, \, \zeta)|^2\, e^{-2m\varphi(z)}\, e^{-2m\varphi(\zeta)} \, d\, \lambda(z)\bigg| \leq $  \\

\hspace{5ex}$\displaystyle \leq \frac{M}{m^n} \, \int\limits_{\Omega} |K_{m\varphi}(z, \, \zeta)|^2\, e^{-2m\varphi(z)}\, e^{-2m\varphi(\zeta)} \, d\, \lambda(z) = \frac{M}{m^n} \, B_m(\zeta)\, e^{-2m\varphi(\zeta)}  $  \\

\hspace{5ex} $\displaystyle \leq 2\, M \, \frac{2^n}{\pi^n} \, \lambda_1(\zeta) \dots \lambda_n(\zeta)$, \hspace{6ex}  $\zeta \in \Omega \setminus V$, \\

\noindent where the equality above follows from the reproducing property of the Bergman kernel, and the last expression is locally integrable on $\Omega \setminus V$, by our key observation. By dominated convergence, it is then enough to prove that:  \\

$\displaystyle\lim\limits_{m\rightarrow +\infty} \frac{1}{m^n} \int\limits_{\Omega} g(z, \, \zeta)\, |K_{m\varphi}(z, \, \zeta)|^2\, e^{-2m\varphi(z)}\, e^{-2m\varphi(\zeta)} \, d\, \lambda(z) = $  \\

\hspace{20ex} $\displaystyle = \frac{2^n}{\pi^n}\, g(\zeta, \, \zeta)\, \lambda_1(\zeta) \dots \lambda_n(\zeta), \hspace{3ex} \zeta \in \Omega \setminus V.$  \\

\noindent To this end, we can repeat Berndtsson's arguments showing the Bergman kernel to decay rapidly off the diagonal, namely that for every $\varepsilon >0$ and every $\zeta \in \Omega\setminus V$, we have: \\

$\displaystyle \lim\limits_{m\rightarrow +\infty}\frac{1}{m^n} \int\limits_{|z-\zeta|>\varepsilon}  |K_{m\varphi}(z, \, \zeta)|^2\, e^{-2m\varphi(z)}\, e^{-2m\varphi(\zeta)} \, d\, \lambda(z)=0$. \\

\noindent The reproducing property of the Bergman kernel then leads to $(b)$. Point $(c)$ is an easy consequence of $(a)$ and $(b)$. \\

 The whole proof of Theorem \ref{The:integercoeff} thus boils down to establishing the uniform upper estimate (\ref{eqn:unifsingestimate}). We will proceed in several steps. \\

\noindent {\bf Step 1.} Assume $\varphi = \psi + \log|g|$ on $\Omega$, for some $g\in {\cal O}(\Omega)$ such that $\mbox{div}\, g$ is a normal crossing divisor, and for some smooth and strictly psh $\psi$ on $\Omega$. Then $V=\{g=0\}$ and $i\partial\bar{\partial}\varphi = i\partial\bar{\partial}\psi$ on $\Omega \setminus V$. \\

 Thus $f\in {\cal O}(\Omega)$ satisfies $1=\int_{\Omega}|f|^2 \, e^{-2m\varphi} = \int_{\Omega} \frac{|f|^2}{|g^m|^2} \, e^{-2m\psi}$ if and only if $f=g^m \, h$ for some $h\in {\cal O}(\Omega)$ satisfying $\int_{\Omega}|h|^2 \, e^{-2m\psi}=1$. This means that $B_{m\varphi} = |g|^{2m} \,  B_{m\psi}$, and we get: \\

\hspace{6ex} $\displaystyle \frac{n!}{m^n}\, B_{m\varphi}\, e^{-2m\varphi} = \frac{n!}{m^n}\, B_{m\psi}\, e^{-2m\psi},$  \hspace{2ex} on $\Omega.$ \\

\noindent This last expression satisfies the uniform upper estimate claimed in (\ref{eqn:unifsingestimate}) on $B\subset\subset \Omega$ thanks to the Berndtsson-Lindholm inequality (\ref{eqn:unifestimate}) applied to the smooth function $\psi$. Notice that, in particular, this yields: \\

\hspace{6ex} $v_B(\psi + \log |g|) = v_B(\psi)$.   \\

\noindent {\bf Step 2.} Assume $\varphi$ has locally divisorial singularities; namely, every point in $\Omega$ has a neighbourhood on which $\varphi = \psi + \log|g|$ for some holomorphic function $g$ such that $\mbox{div}\, g$ has normal crossings, and for some smooth strictly psh function $\psi$. \\

 Let $U\subset\subset \Omega$ be a pseudoconvex such neighbourhood. It is clear that the unit ball of ${\cal H}_{\Omega}(m\varphi)$ injects into the unit ball of ${\cal H}_U(m\varphi_{|U})$ under restriction to $U$. Therefore, in the light of (\ref{eqn:defBergman}), we get: \\

$\displaystyle \frac{n!}{m^n}\, B_{m\varphi}\, e^{-2m\varphi} \leq \frac{n!}{m^n}\, B_{m\varphi_{|U}}\, e^{-2m\varphi}= \frac{n!}{m^n}\, B_{m\psi_{|U}}\, e^{-2m\psi},$  \hspace{3ex} on $U.$  \\

\noindent  The last term satisfies the uniform upper estimate claimed in (\ref{eqn:unifsingestimate}) on $U\setminus V$ thanks again to the smooth case applied to $\psi$. The resulting constant $C>0$ depends on $U$, but we get a constant independent of $m$ for the estimate on $B\setminus V$ after taking a finite covering of $\bar{B}$ by such open sets $U$.  \\

\noindent {\bf Step 3.} Assume $\varphi = \frac{p}{2}\, \log(|g_1|^2 + \dots + |g_N|^2) + u$ on $\Omega$, for some holomorphic functions $g_1, \dots , \, g_N$ and some smooth function $u$. \\

Let $J:=(g_1, \dots , \, g_N) \subset {\cal O}_{\Omega}$ be the ideal sheaf generated by $g_1, \dots , \, g_N$, and let $\mu : \tilde{\Omega} \rightarrow \Omega$ be a proper modification such that \\

\hspace{6ex} $\mu^{\star}\, J = {\cal O}(-E)$, \\

\noindent $\tilde{\Omega}$ is a smooth manifold, and $E$ is an effective normal crossing divisor on $\tilde{\Omega}$. The divisor $E$ can be locally represented as $E=\mbox{div}\, g$, for some locally defined holomorphic function $g$. We then get, locally on $\tilde{\Omega}$: \\

$\displaystyle \varphi\circ \mu = \frac{p}{2} \, \log\bigg(\bigg|\frac{g_1\circ \mu}{g}\bigg|^2 + \dots + \bigg|\frac{g_N\circ \mu}{g}\bigg|^2\bigg) + u\circ \mu + \log |g^p|$ \\

\hspace{5ex} $= \psi + \log |g^p|,$ \\

\noindent where $\psi$ is a smooth strictly psh function. Thus $\varphi\circ \mu$ has locally divisorial singularities and satisfies the hypotheses of Step 2. In passing from Bergman kernels defined on bigger sets to Bergman kernels defined on smaller sets, we shall need the following comparison lemma.

\begin{Lem}\label{Lem:comparison} If $U\subset\subset \Omega$ is a pseudoconvex open subset, the Bergman kernels associated with the weight $m\varphi$ on $\Omega$ and respectively $U$:  \\

\hspace{6ex} $\displaystyle B_{m\varphi, \, \Omega}:=\sum\limits_{k=0}^{+\infty}|\sigma_{m, \, k}|^2,$ \hspace{3ex} $\displaystyle B_{m\varphi, \, U}:=\sum\limits_{k=0}^{+\infty}|\mu_{m, \, k}|^2,$ \\

\noindent defined by orthonormal bases $(\sigma_{m, \, k})_{k\in \N}$ and $(\mu_{m, \, k})_{k\in \N}$ of the Hilbert spaces ${\cal H}_{\Omega}(m\varphi)$ and respectively ${\cal H}_{U}(m\varphi)$, can be compared, for every $m\in \N$, as: \\

\hspace{3ex} $\displaystyle B_{m\varphi, \, \Omega} \leq B_{m\varphi, \, B} \leq C_{n, \, d, \, r}\, B_{m\varphi, \, \Omega}$ \hspace{2ex} on any $U_0\subset\subset U\subset\subset \Omega,$ \\

\noindent where $C_{n, \, d, \, r}>0$ is a constant depending only on $n$, the diameter $d$ of $\Omega$, and the distance $r>0$ between the boundaries of $U_0$ and $U$.

\end{Lem}

\noindent {\it Proof.} As the restriction to $U$ defines an injection of unit ball of ${\cal H}_{\Omega}(m\varphi)$ into the unit ball of ${\cal H}_{B}(m\varphi)$, the former inequality follows. For the latter inequality, fix $x\in U_0$ and let $f\in {\cal O}(U)$ such that $\int_U |f|^2\, e^{-2m\varphi}=1$ be an arbitrary element in the unit sphere of ${\cal H}_U(m\varphi)$. We use H\"ormander's $L^2$ estimates ([Hor65]) to construct a holomorphic function $F\in {\cal H}_{\Omega}(m\varphi)$ such that $F(x)=f(x)$ and \\

\hspace{6ex} $\displaystyle \int\limits_{\Omega}|F|^2\, e^{-2m\varphi} \leq 2\bigg(1+C_n\, \frac{d^{2n}\, e^{2d^2}}{r^{2(n+1)}}\bigg)\, \int\limits_U |f|^2\, e^{-2m\varphi}=C_{n, \, d,\, r}.$ \\

\noindent This is done by solving the equation: \\

\hspace{6ex}$\bar{\partial}u = \bar{\partial}(\theta\, f) \hspace{2ex} \mbox{on} \,\, \Omega$

\vspace{2ex}

 \noindent with a cut-off function $\theta\in C^{\infty}(\C^n)$, $\mbox{Supp}\, \theta \subset\subset U$, $\theta\equiv 1$ on $U_0\subset\subset U$, and with the strictly psh weight $m\varphi + n\log|z-x| + |z-x|^2$. There exists a $C^{\infty}$ solution $u$ satisfying the estimate: \\

\hspace{6ex} $\displaystyle \int\limits_{\Omega}\frac{|u|^2}{|z-x|^{2n}}\, e^{-2m\varphi}\, e^{-2|z-x|^2} \leq 2 \int\limits_{\Omega}\frac{|\bar{\partial}\theta|^2\, |f|^2}{|z-x|^{2n}}\, e^{-2m\varphi}\, e^{-2|z-x|^2}.$ \\

\noindent Due to the non-integrability of $|z-x|^{-2n}$ near $x$, we have $u(x)=0$. Thus $F$ is obtained as:\\

\hspace{6ex} $F:=\theta \, f -u\in {\cal O}(\Omega).$ \\

\noindent Now $F/\sqrt{C_{n, \, d, \, r}}$ belongs to the unit ball of ${\cal H}_{\Omega}(m\varphi)$, and we get: \\

\hspace{6ex} $\displaystyle |f(x)|^2=|F(x)|^2 \leq C_{n, \, d, \, r}\, B_{m\varphi, \, \Omega}(x),$ \\

\noindent which proves the latter inequality by taking the supremum over all $f$ in the unit sphere of ${\cal H}_U(m\varphi)$. \hfill $\Box$

\vspace{2ex}

\noindent We now resume Step 3 of the proof of Theorem \ref{The:integercoeff}. We may assume, without loss of generality, that the Jacobian $J_{\mu}$ of $\mu$ is globally defined on $\tilde{\Omega}$. Otherwise we can work on smaller open subsets of $\tilde{\Omega}$ contained in coordinate patches and the previous Lemma \ref{Lem:comparison} shows that the Bergman kernels are only distorted by an insignificant constant $C_{n,\,d,\,r}$ independent of $m$. A change of variable shows that, for each $\sigma_{m, \, j}$, we have: \\

$\displaystyle 1 = \int\limits_{\Omega} |\sigma_{m, \, j}|^2 \, e^{-2m\varphi}\, d\lambda = \int\limits_{\Omega\setminus V} |\sigma_{m, \, j}|^2 \, e^{-2m\varphi}\, d\, \lambda$  \\

\hspace{2ex} $\displaystyle  = \int\limits_{\tilde{\Omega}\setminus \mbox{Supp}\, E} |\sigma_{m, \, j}\circ \mu|^2 \, |J_{\mu}|^2 \, e^{-2m\varphi\circ \mu}\, d\tilde{V} = \int\limits_{\tilde{\Omega}} |\sigma_{m, \, j}\circ \mu|^2 \, |J_{\mu}|^2 \, e^{-2m\varphi\circ \mu}\, d\tilde{V},$\\

\noindent for a suitable volume form $d\, \tilde{V}$ on $\tilde{\Omega}$. Consequently, if $\displaystyle B_{m\varphi\circ \mu}$ is the Bergman kernel associated with $m\, \varphi\circ \mu$ and the volume form $d\, \tilde{V}$ on $\tilde{\Omega}$, we have \\

 \hspace{10ex} $\displaystyle B_{m\varphi\circ \mu}=|J_{\mu}|^2 \, B_{m\varphi}\circ \mu$.\\

\noindent  Thus, proving the upper estimate claimed in (\ref{eqn:unifsingestimate}) amounts to proving:

\begin{equation}\label{eqn:comparison}\displaystyle \frac{n!}{m^n} \, B_{m \varphi\circ \mu}\, e^{-2m\varphi\circ \mu} \leq C\, \frac{2^n}{\pi^n} \, |J_{\mu}|^2 \, \lambda_1 \circ \mu \dots \lambda_n \circ \mu,\end{equation}

\vspace{2ex}

\noindent on $\mu^{-1}(B)\setminus \mbox{Supp}\, E$, since $\mu$ defines an isomorphism between $\tilde{\Omega}\setminus \mbox{Supp}\, E$ and $\Omega \setminus V$, and $J_{\mu}$ does not vanish on $\tilde{\Omega} \setminus \mbox{Supp}\, E$. If $\tilde{\lambda_1}, \dots , \, \tilde{\lambda_n}$ are the eigenvalues of $i\partial\bar{\partial}(\varphi \circ \mu)$, it can easily be seen that \\

\hspace{10ex}  $|J_{\mu}|^2 \, \lambda_1 \circ \mu \dots \lambda_n \circ \mu = \tilde{\lambda_1} \dots \tilde{\lambda_n}.$ \\

\noindent Meanwhile, $\varphi\circ \mu$ has locally divisorial singularities, and Step 2 can be applied to get the uniform estimate (\ref{eqn:comparison}) which proves the uniform estimate claimed in (\ref{eqn:unifsingestimate}) and completes the proof of Theorem \ref{The:integercoeff}.  \hfill $\Box$

\vspace{3ex}

\noindent {\bf $(b)$ \, \, Analytic singularities with arbitrary coefficients} \\

Let us first consider the case of divisorial singularities with noninteger coefficients. For the sake of simplicity, we assume that our domains are polydiscs, $\Omega=D^n$ and $B=D_{1-\varepsilon}^n$, where $D$ (respectively $D_{1-\varepsilon}$) is the unit disc (respectively the disc of radius $1-\varepsilon$) in $\C$.

\begin{Prop}\label{Prop:arbitrarycoeff} Assume that $\varphi(z) = \psi(z) + c_1 \, \log|z_1| + \dots + c_n \, \log|z_n|$, 

\noindent with $\psi(z)=\psi_1(z_1) + \dots + \psi_n(z_n)$ for some $C^{\infty}$ functions $\psi_j$ on $\C$ depending only on $|z_j|$ respectively, and some constants $c_j>0$, $j=1, \dots , \, n$. If $i\partial\bar{\partial}\varphi \geq C_0 \, \omega$ for some $C_0>0$, then: \\

\hspace{6ex} $v_B(\psi + \sum\limits_{j=1}^{n} c_j \, \log|z_j|)= v_B(\psi)=\displaystyle \frac{2^n}{\pi^n}\, \int_B (i\partial\bar{\partial}\psi)^n <+\infty.$

\end{Prop}

\noindent {\it Proof.} Only the first equality needs a proof. The second equality follows from Theorem \ref{The:Bo} for smooth functions. To begin with, we shall prove that:

\begin{equation}\label{eqn:integralcomparison}\displaystyle \int\limits_B B_{m\varphi}\, e^{-2m\varphi}\, d\,\lambda \geq  \int\limits_B B_{m\psi}\, e^{-2m\psi}\, d\,\lambda, \hspace{2ex} \mbox{for every} \hspace{2ex} m\in \N^{\star},\end{equation}

\noindent which clearly implies that

\begin{equation}\label{eqn:increase} v_B(\psi + \sum\limits_{j=1}^{n} c_j \, \log|z_j|)\geq v_B(\psi).\end{equation}

\noindent Fix $m\in \N^{\star}$, and let $a_j:=\{m\, c_j\}$ for $j=1, \dots , \, n$, where $\{ \,\cdot \, \}$ stands for the fractional part. As already noticed at Step 1 in the proof of Theorem \ref{The:integercoeff}, we have: \\

$B_{m\varphi}\, e^{-2m\varphi} = B_{m\psi + a_1\, \log |z_1| + \dots + a_n \log|z_n|}\, e^{-2m\psi - a_1\, \log |z_1|- \dots - a_n\, \log |z_n|}$. \\

\noindent This already implies that $v_B(\psi + \sum\limits_{j=1}^{n} c_j \, \log|z_j|)=v_B(\psi + \sum\limits_{j=1}^{n} \{c_j\} \, \log|z_j|)$, and thus we may assume, without loss of generality, that $0\leq c_j <1$ for $j=1, \dots , \, n$.

\noindent Since $a_j<1$ for every $j$, the exponential $e^{-2m\varphi}$ is easily seen to be locally integrable (the integral being a product of integrals depending each on one complex variable). Thus all the monomials $z^{\alpha}=z_1^{\alpha_1} \dots z_n^{\alpha_n}$, $\alpha=(\alpha_1, \dots , \alpha_n)\in \N^n,$ make up a complete orthogonal set of ${\cal H}_{\Omega}(m\varphi)$, which becomes an orthonormal basis after each monomial is normalized to have norm $1$. We thus get: \\

$\displaystyle B_{m\varphi}\, e^{-2m\varphi} = \sum\limits_{\alpha_1, \dots, \, \alpha_n=0}^{+\infty} \prod\limits_{j=1}^n \frac{|z_j|^{2(\alpha_j - a_j)}\, e^{-2m\psi_j}}{\displaystyle \int_D |z_j|^{2(\alpha_j - a_j)}\, e^{-2m\psi_j(z_j)}\, d\, \lambda_1(z_j)},$ \\

\noindent and the analogous formula for $m\psi$: \\

$\displaystyle B_{m\psi}\, e^{-2m\psi} = \sum\limits_{\alpha_1, \dots, \, \alpha_n=0}^{+\infty} \prod\limits_{j=1}^n \frac{|z_j|^{2\alpha_j}\, e^{-2m\psi_j}}{\displaystyle \int_D |z_j|^{2\alpha_j}\, e^{-2m\psi_j(z_j)}\, d\, \lambda_1(z_j)}.$ \\

\noindent It is then enough to prove that, for every $j$, we have: 

\begin{equation}\label{eqn:fractioncomparison} \displaystyle \frac{\displaystyle \int_{D_{1-\varepsilon}}|z_j|^{2(\alpha_j - a_j)}\, e^{-2m\psi_j}\, d\,\lambda_1(z_j)}{\displaystyle \int_D |z_j|^{2(\alpha_j - a_j)}\, e^{-2m\psi_j}\, d\,\lambda_1(z_j)} \geq \frac{\displaystyle \int_{D_{1-\varepsilon}}|z_j|^{2\alpha_j}\, e^{-2m\psi_j}\, d\,\lambda_1(z_j)}{\displaystyle \int_D |z_j|^{2\alpha_j}\, e^{-2m\psi_j}\, d\,\lambda_1(z_j)}.\end{equation}

\noindent Taking polar coordinates $z_j=r_j \, \rho_j$, $0\leq r_j \leq 1$, $\rho_j\in S^1$, we are then reduced to proving that, given a $C^{\infty}$ function $u\geq 0$ and a constant $0\leq c<1$, we have: \\

$\frac{\displaystyle \int_0^{1-\varepsilon} x^{2k+1} \, u(x)\, \frac{1}{x^{2c}}\, dx}{\displaystyle \int_0^1 x^{2k+1} \, u(x)\, \frac{1}{x^{2c}}\, dx} \geq \frac{\displaystyle \int_0^{1-\varepsilon} x^{2k+1} \, u(x)\, dx}{\displaystyle \int_0^1 x^{2k+1} \, u(x)\, dx},$ \\

\noindent or equivalently, that: \\

$\frac{\displaystyle \int_0^{1-\varepsilon} x^{2k+1} \, u(x)\, \frac{1}{x^{2c}}\, dx}{\displaystyle \int_0^{1-\varepsilon} x^{2k+1} \, u(x)\, dx} \geq \frac{\displaystyle \int_{1-\varepsilon}^1 x^{2k+1} \, u(x)\, \frac{1}{x^{2c}}\, dx}{\displaystyle \int_{1-\varepsilon}^1 x^{2k+1} \, u(x)\, dx},$ \\

\noindent which is clear since the left-hand side is $\geq \frac{1}{(1-\varepsilon)^{2c}}$, and the right-hand side is $\leq \frac{1}{(1-\varepsilon)^{2c}}.$ Inequalities (\ref{eqn:integralcomparison}) and (\ref{eqn:increase}) are thus proved. \\

 We shall now prove that:

\begin{equation}\label{eqn:increasebis} v_B(\psi + \sum\limits_{j=1}^{n} c_j \, \log|z_j|)\leq v_B(\psi + \sum\limits_{j=1}^{n} \log|z_j|).\end{equation}

\noindent Let $v(z):= \psi(z) + \log|z_1| + \dots + \log|z_n|$. The monomials $z^{\alpha}=z_1^{\alpha_1}\dots z_n^{\alpha_n}$, with $\alpha_1, \dots , \alpha_n \geq m$, make up a complete orthogonal system of ${\cal H}_{\Omega}(mv)$. Therefore, \\

$\displaystyle B_{mv}(z)\, e^{-2mv(z)} = \sum\limits_{\alpha_1, \dots , \alpha_n =1}^{+\infty} \frac{\displaystyle |z_1|^{2(\alpha_1-1)} \dots |z_n|^{2(\alpha_n-1)}}{\displaystyle \int_{\Omega} |z_1|^{2(\alpha_1-1)} \dots |z_n|^{2(\alpha_n-1)}\, e^{-2m\psi(z)}} \, e^{-2m\psi(z)} $ \\

\hspace{13ex} $\displaystyle = \sum\limits_{\alpha_1, \dots , \alpha_n =1}^{+\infty} \prod\limits_{j=1}^n \frac{\displaystyle |z_j|^{2(\alpha_j-1)}\, e^{-2m\psi_j(z_j)}}{\displaystyle \int_D |z_j|^{2(\alpha_j-1)}\, e^{-2m\psi_j(z_j)}\, d\,\lambda_1(z_j)},$ \\

\noindent which entails: \\

$\displaystyle \frac{1}{m^n} \, \int_B B_{mv}(z)\, e^{-2mv(z)} \, d\,\lambda(z) = \frac{1}{m^n} \prod\limits_{j=1}^n \sum\limits_{\alpha_j=1}^{+\infty} \frac{\displaystyle \int_{D_{1-\varepsilon}} |z_j|^{2(\alpha_j-1)}\, e^{-2m\psi_j(z_j)}\, d\, \lambda_1(z_j)}{\displaystyle \int_D |z_j|^{2(\alpha_j-1)}\, e^{-2m\psi_j(z_j)}\, d\,\lambda_1(z_j)}$ \\

$\displaystyle  \geq \frac{1}{m^n} \prod\limits_{j=1}^n \bigg(\sum\limits_{\alpha_j=0}^{+\infty} \frac{\displaystyle \int_{D_{1-\varepsilon}} |z_j|^{2(\alpha_j-c_j)}\, e^{-2m\psi_j(z_j)}\, d\, \lambda_1(z_j)}{\displaystyle \int_D |z_j|^{2(\alpha_j-c_j)}\, e^{-2m\psi_j(z_j)}\, d\,\lambda_1(z_j)} - \frac{\displaystyle \int_{D_{1-\varepsilon}} |z_j|^{-2c_j}\, e^{-2m\psi_j(z_j)}\, d\, \lambda_1(z_j)}{\displaystyle \int_D |z_j|^{-2c_j}\, e^{-2m\psi_j(z_j)}\, d\,\lambda_1(z_j)} \bigg)$ \\

$\displaystyle \geq \frac{1}{m^n} \prod\limits_{j=1}^n \bigg(\sum\limits_{\alpha_j=0}^{+\infty} \frac{\displaystyle \int_{D_{1-\varepsilon}} |z_j|^{2(\alpha_j-c_j)}\, e^{-2m\psi_j(z_j)}\, d\,\lambda_1(z_j)}{\displaystyle \int_D |z_j|^{2(\alpha_j-c_j)}\, e^{-2m\psi_j(z_j)}\, d\,\lambda_1(z_j)} - 1 \bigg)$ \\

$\displaystyle \geq  \frac{1}{m^n}\,\int_B B_{m\varphi}\, e^{-2m\varphi}\, d\,\lambda_n - \frac{1}{m}\, \sum\limits_{k=1}^n \frac{1}{m^{n-1}}\, \int_{D_{1-\varepsilon}^{n-1}} B_{mu_k}\, e^{-2mu_k}\, d\,\lambda_{n-1}$ \\

 $\displaystyle + \frac{1}{m^2}\, \sum\limits_{k_1, \, k_2} \frac{1}{m^{n-2}}\, \int_{D_{1-\varepsilon}^{n-2}} B_{mu_{k_1, \, k_2}}\, e^{-2mu_{k_, \, k_2}}\, d\,\lambda_{n-2} -  \dots + \frac{(-1)^n}{m^n},$  \\

\noindent where we have denoted $u_k(z_1, \dots , \hat{z_k}, \dots z_n):= \varphi(z) -\psi_k(z_k) - c_k\, \log|z_k|$, and the analogous expressions when several indices are missing. The $k$-dimensional Lebesgue measure has been denoted $d\, \lambda_k$. Note that as we assume $0\leq c_j <1$, the first inequality above is implied by estimate (\ref{eqn:fractioncomparison}) with $\alpha_j$ replaced with $\alpha_j - c_j$, and $c_j$ replaced with $1-c_j$.

\vspace{3ex}

\noindent We can thus run an induction on the dimension $n$. If we assume the finiteness of the volume for psh functions of the form under consideration defined in less than $n$ variables, all the terms appearing on the right-hand side, except the first one, tend to $0$ as $m\rightarrow +\infty$. Inequality (\ref{eqn:increasebis}) is then what we get in the limit. \\

\noindent  Now, as already noticed (e.g. at Step 1 above), \\

 \hspace{10ex} $v_B(\psi) = v_B(\psi + \sum\limits_{j=1}^n\log|z_j|)$,  \\

\noindent which, alongside inequalities (\ref{eqn:increase}) and (\ref{eqn:increasebis}), completes the proof.  \hfill $\Box$

\vspace{2ex}

 We can now round off the study of the finiteness of the volume of a psh function in the case of general analytic singularities with arbitrary (not necessarily integer) coefficients.

\begin{Prop}\label{Prop:general} Let $\varphi = \frac{c}{2} \, \log(|g_1|^2 + \dots + |g_N|^2) + u$, for some holomorphic functions  $g_1, \dots , \, g_N$, some $c>0$, and some $C^{\infty}$ function $u$ on $\Omega\subset\subset \C^n$. Assume $i\partial\bar{\partial}\varphi \geq C_0 \, \omega$ for some constant $C_0>0$. Then: \\

\hspace{6ex} $\displaystyle v_B(\varphi):=\limsup\limits_{m\rightarrow +\infty} \frac{n!}{m^n}\, \int\limits_{B}B_{m\varphi}\, e^{-2m\varphi}\, d\, V <+\infty.$  

\end{Prop}

\noindent {\it Proof.} Let $J=(g_1, \dots , \, g_N)\subset {\cal O}_{\Omega}$, and let $\mu : \tilde{\Omega} \rightarrow \Omega$ be a proper modification such that $\tilde{\Omega}$ is smooth and $\mu^{\star}\, J = {\cal O}(-E)$ for an effective normal crossing divisor $E$ on $\tilde{\Omega}$. The change of variable formula shows, as at Step 3 in the proof of Theorem \ref{The:integercoeff}, that \\

\hspace{6ex} $v_B(\varphi)=v_{\mu^{-1}(B)}(\varphi\circ \mu)$. \\

\noindent  If we cover $\overline{\mu^{-1}(B)}$ by finitely many open polydiscs $U_k$ such that $E=\mbox{div} g_k$ on $U_k$, and use the observation made at Step 2 in the proof of Theorem \ref{The:integercoeff} on the behaviour of Bergman kernels under restrictions, we see that: \\

\hspace{6ex} $\displaystyle v_{\mu^{-1}(B)}(\varphi\circ \mu) \leq \sum\limits_{k} v_{U_k}(\varphi \circ \mu_{|U_k}).$ \\

\noindent Now, $\varphi \circ \mu_{|U_k}$ can be written in the form of functions $\varphi$ considered in the previous Proposition \ref{Prop:arbitrarycoeff}. Therefore, $v_{U_k}(\varphi \circ \mu_{|U_k})< +\infty$ for every $k$. This completes the proof. \hfill $\Box$

\vspace{2ex}

 The conclusion of these considerations about the volume of a psh function with analytic singularities is that if $\lambda_{m, \, 1}, \dots ,\, \lambda_{m, \, N_m}$ are those eigenvalues of  the concentration operator $T_{B, \, m}$ that are $\geq 1-\varepsilon$, for some $0<\varepsilon <1$, then their number can be estimated as: \\

\begin{equation}\label{eqn:numberestimate}\displaystyle N_m \leq \frac{v_B(\varphi)}{1-\varepsilon}\, \frac{m^n}{n!} = O(m^n), \hspace{3ex} \mbox{for} \hspace{2ex} m>>1.\end{equation}

\subsection {Effective local finite generation} \label{subsection:generation}

 In this section, we will prove Theorem \ref{The:introd1}. Let $(\sigma_{m, \, j})_{j\in \mathbb{N}^{\star}}$ be an orthonormal basis of ${\cal H}_{\Omega}(m\varphi)$ made up of eigenvectors of the Toeplitz concentration operator $T_{B, \, m}$ discussed in the previous section. Recall that we are looking for finitely many among these that generate the ideal sheaf ${\cal I}(m\varphi)$ on $B\subset\subset \Omega$.

\begin{Def}\label{Def:concentration} A nonzero function $f\in {\cal H}_{\Omega}(m\varphi)$ is said to be $\delta$-concentrated on $B$, for some $\delta >0$, if $\displaystyle \frac{\int_B|f|^2 \, e^{-2m\varphi}}{\int_{\Omega}|f|^2 \, e^{-2m\varphi}} \geq 1-\delta$. 

\end{Def}

Since the restrictions to $B$ of the $\sigma_{m, \, j}$'s are still orthogonal to one another (cf. lemma \ref{Lem:concentratedbasis}), an element $f$ in the subspace of ${\cal H}_{\Omega}(m\varphi)$ generated by $\sigma_{m, \, 1}, \dots , \sigma_{m, \, N_m}$ is $\varepsilon$-concentrated on $B$, like its generators. The following simple observation shows that the converse is not far from being true.

\begin{Lem}\label{Lem:concentration} Let $f\in {\cal H}_{\Omega}(m\varphi)$ be $\varepsilon^2$-concentrated on $B$. Then, there is an element $g$ in the subspace of ${\cal H}_{\Omega}(m\varphi)$ generated by the $\varepsilon$-concentrated $\sigma_{m, \, 1}, \dots , \sigma_{m, \, N_m}$ such that $\int_{\Omega}|f-g|^2 \, e^{-2m\varphi} < \varepsilon \, \int_{\Omega}|f|^2 \, e^{-2m\varphi}$.

\end{Lem}

\noindent {\it Proof.} We may assume that $\int_{\Omega}|f|^2 \, e^{-2m\varphi} = 1$. Let $f= \sum\limits_{j=1}^{N_m}a_j \, \sigma_{m, \, j} + \sum\limits_{k=N_m +1}^{+\infty}a_k \, \sigma_{m, \, k}$ be the decomposition of $f$ with respect to the chosen orthonormal basis of ${\cal H}_{\Omega}(m\varphi)$, where $a_j, \, a_k \in \mathbb{\C}$. Then, $\sum\limits_{j=1}^{+\infty}|a_j|^2 =1$, and \\

$\displaystyle \int_B|f|^2 \, e^{-2m\varphi} = \sum\limits_{j=1}^{+\infty}|a_j|^2 \, \int_B|\sigma_{m, \, j}|^2 \, e^{-2m\varphi}$. We then get: \\

$\displaystyle \int_{\Omega\setminus B}|f|^2 \, e^{-2m\varphi} = \sum\limits_{j=1}^{+\infty}|a_j|^2 \, \int_{\Omega\setminus B}|\sigma_{m, \, j}|^2 \, e^{-2m\varphi}  \leq \varepsilon^2$.  \\

\noindent Since $\int_{\Omega\setminus B}|\sigma_{m, \, k}|^2\, e^{-2m\varphi} > \varepsilon$, for every $k\geq N_m+1$, we get $\sum\limits_{k=N_m +1}^{+\infty}|a_k|^2 < \varepsilon$. If we set $g:=\sum\limits_{j=1}^{N_m}a_j \, \sigma_{m, \, j}$, the lemma is proved.  \hfill $\Box$

\vspace{3ex}

This strongly suggests where to turn for the most likely choice of finitely many local generators for ${\cal I}(m\varphi)$. For a small given $\varepsilon>0,$ let $\sigma_{m, \, 1}, \dots , \, \sigma_{m, \, N_m}$ be those elements which are $\varepsilon$-concentrated on $B$ (cf. Definition \ref{Def:concentration}). We will now show that, if $\varepsilon$ is well chosen and $m$ is big enough, the ideal sheaf ${\cal I}(m\varphi)$ is generated, on a relatively compact open subset, by $\sigma_{m, \, 1}, \dots , \, \sigma_{m, \, N_m},$ with an effective control of the coefficients. We shall proceed inductively. The crux is the following approximation to order one of a local section of ${\cal I}(m\varphi)$ by a finite linear combination of $\sigma_{m, \, j}$'s. The standard K\"ahler form on $\C^n$ is denoted throughout by $\omega$.

 \begin{Lem}\label{Lem:firstordergeneration} Let $\varphi$ be a strictly psh function on a bounded pseudoconvex open set $\Omega\subset \C^n$ such that $i\partial\bar{\partial}\varphi \geq C_0\, \omega$ for some constant $C_0>0$. Let $B:= B(x, \, r)\subset\subset\Omega$ be an arbitrary open ball. Then, there exist $\varepsilon>0$ and $m_0=m_0(C_0)\in \mathbb{N}$ such that, for every $m\geq m_0,$ the following property holds. Every $g\in {\cal H}_B(m\varphi)$ admits, with respect to the $\varepsilon$-concentrated $\sigma_{m, \, j}$'s, a decomposition:  \\

$\displaystyle g(z) = \sum\limits_{j=1}^{N_m} c_j\, \sigma_{m, \, j}(z) + \sum\limits_{l=1}^{n} (z_l-x_l) \, h_l(z) $,  \hspace{3ex}  $z\in B,$  \\

\noindent with some $c_j\in \C$ satisfying $\displaystyle \sum\limits_{j=1}^{N_m}|c_j|^2 \leq C\, N_m\, \int_B|g|^2\, e^{-2m\varphi}$, and some holomorphic functions $h_l$ on $B$, satisfying: \\

$\displaystyle \sum\limits_{l=1}^{n} \int_{B} |h_l|^2 \, e^{-2m\varphi} \leq C \, \int_B|g|^2\, e^{-2m\varphi},$  \\

\noindent where $C>0$ is a constant depending only on $n,$ on $r,$ and on the diameter of $\Omega.$ 

\end{Lem}

\noindent {\it Proof.} We may assume, without loss of generality, that $x=0.$ Let $\theta\in C^{\infty}(\Omega), 0\leq \theta \leq 1,$ be a cut-off function such that $\mbox{Supp} \, \theta \subset B$, $\theta\equiv 1$ on some arbitrary $B':=B(0, \, r') \subset\subset B$, and $|\bar{\partial}\theta|\leq \frac{3}{r}.$ Let us fix some $g\in {\cal H}_B(m\varphi)$ such that $C_g:=\int_B|g|^2\, e^{-2m\varphi}<+\infty$. We will now use H\"ormander's $L^2$ estimates ([Hor65]) to solve the equation:  

\vspace{1ex}

\hspace{6ex} $\bar{\partial}u = \bar{\partial}(\theta g),$  \hspace{3ex} on $\Omega,$ 

\vspace{1ex}

\noindent with the weight $m\varphi(z) + (n+1) \, \log|z|$. We get a solution $u\in C^{\infty}(\Omega)$ such that: 

\begin{equation}\label{eqn:uestimate}\displaystyle \int_{\Omega}\frac{|u|^2}{|z|^{2(n+1)}}\, e^{-2m\varphi} \leq \frac{1}{C_0\,m} \, \int_{\Omega}\frac{|\bar{\partial}\theta|^2 \, |g|^2}{|z|^{2(n+1)}} \, e^{-2m\varphi}\leq \frac{9}{r^2}\, \frac{1}{C_0\,m} \, \int_{B\setminus B'}\frac{|g|^2}{|z|^{2(n+1)}} \, e^{-2m\varphi}.\end{equation}

\noindent Put $F_m:=\theta g - u\in{\cal H}_{\Omega}(m\varphi)$ and get the decomposition $g= F_m + (g-F_m)$ on $B$. The above estimate for $u$ implies:  

\begin{equation}\label{eqn:Skodaestimate}\displaystyle \int_B \frac{|g-F_m|^2}{|z|^{2(n+1)}} \, e^{-2m\varphi} \leq 2 \, \int_{B\setminus B'} \frac{|1-\theta|^2 |g|^2}{|z|^{2(n+1)}} \, e^{-2m\varphi} + 2 \, \int_B\frac{|u|^2}{|z|^{2(n+1)}}\, e^{-2m\varphi}\end{equation}

\hspace{19ex} $\leq C(r, \, r', \, d)\, C_g$,  \\

\noindent with a constant $C(r, \, r', \, d)>0$ depending only on $r, \, r',$ and the diameter $d$ of $\Omega$, if $m$ is chosen so big that $\frac{1}{C_0\, m}<1.$

\noindent We shall now apply Skoda's $L^2$ division theorem (cf. [Sko72]) to obtain: 

\begin{equation}\label{eqn:Skoda}\displaystyle g(z)-F_m(z) = \sum\limits_{l=1}^{n}z_l \, v_l(z),  \hspace{3ex}  z\in B,\end{equation}  

\noindent for some holomorphic functions $v_l$ on $B$ satisfying $\sum_{l=1}^{n} \int_B \frac{|v_l|^2}{|z|^{2n}} \, e^{-2m\varphi} \leq 2\, C(r, \, r', \, d)\, C_g.$ Since $|z|^{2n} \leq r^{2n}$, for $z\in B,$ we get:

\begin{equation}\label{eqn:vs}\displaystyle \sum_{l=1}^{n} \int_B |v_l|^2 \, e^{-2m\varphi} \leq C_1(r, \, r', \, d)\, C_g,\end{equation}

\noindent where $C_1(r, \, r', \, d)=2\, r^{2n}\, C(r, \, r', \, d).$ As for $F_m$, the obvious pointwise inequality $|F_m|^2 \leq 2(|\theta \, g|^2 + |u|^2)$, combined with estimate (\ref{eqn:uestimate}) for $u$, gives: \\

 $\displaystyle \int_B|F_m|^2 \, e^{-2m\varphi} \leq C(r, \, r', \, d) \, C_g$, \\

\noindent after absorbing an extra $d^{2(n+1)}$ in the constant $C(r, \, r', \, d)$. On the other hand, the factor $\frac{1}{C_0\, m}$ in the estimate (\ref{eqn:uestimate}) for $u$ shows that if $m$ is chosen big enough, the $L^2$ norm of $u$ on $\Omega$ is very small in comparison with the $L^2$ norm of $g$ on $B$. This is where the strict psh assumption on $\varphi$ comes in. Since $F_m=g-u$ on $B'$ and $F_m=-u$ on $\Omega \setminus B$, we get, for $m$ big enough and some constant $C_1=C_1(C_0) >0$ independent of $m$:  \\

$\displaystyle 1- \frac{C_1}{m} \leq \frac{\int_B|F_m|^2 \, e^{-2m\varphi}}{\int_{\Omega}|F_m|^2 \, e^{-2m\varphi}} < 1$.  \\

\noindent In other words, $F_m$ is $\frac{C_1}{m}$-concentrated on $B$. Fix some small $\varepsilon>0$ whose choice will be explained later. If $m_{\varepsilon}$ is chosen such that $\frac{C_1}{m}<\varepsilon^2$ for $m\geq m_{\varepsilon},$ then $F_m$ is, in particular, $\varepsilon^2$-concentrated on $B$. Then, lemma \ref{Lem:concentration} shows that in the decomposition:  \\

$\displaystyle F_m(z) = \sum\limits_{j=1}^{N_m} a_j \, \sigma_{m, \, j}(z) + \sum\limits_{k=N_m + 1}^{+\infty}a_k \, \sigma_{m, \, k}(z)$,  \hspace{3ex} $z\in \Omega$, \\

\noindent of $F_m$ with respect to the chosen orthonormal basis of ${\cal H}_{\Omega}(m\varphi)$, we have:

\begin{equation}\label{eqn:sigmas}\displaystyle \sum\limits_{k=N_m + 1}^{+\infty}|a_k|^2 \leq C(r, \, r', \, d)\, \varepsilon\, C_g, \hspace{3ex}  \displaystyle \sum\limits_{j=1}^{N_m}|a_j|^2 \leq  C(r, \, r', \, d) \, C_g.\end{equation}

\noindent If we set $g_1(z):=\sum\limits_{k=N_m + 1}^{+\infty}a_k \, \sigma_{m, \, k}(z)$, (\ref{eqn:Skoda}) gives the decomposition:  

\begin{equation}\label{eqn:g1}\displaystyle g(z) = \sum\limits_{j=1}^{N_m} a_j \, \sigma_{m, \, j}(z) + \sum\limits_{l=1}^{n}z_l \, v_l(z) + g_1(z),  \hspace{3ex}  z\in B,\end{equation} 

\noindent with an effective control of the $a_j$'s and $v_l$'s (see (\ref{eqn:sigmas}) and (\ref{eqn:vs})), and such that: \\

$\displaystyle \int_B |g_1|^2 \, e^{-2m\varphi} \leq \int_{\Omega} |g_1|^2 \, e^{-2m\varphi}= \sum\limits_{k=N_m + 1}^{+\infty}|a_k|^2  \leq C(r, \, r', \, d)\, \varepsilon\, C_g.$  \\

\noindent  Our aim is clearly to get rid of this ``small'' $g_1$ in (\ref{eqn:g1}). The idea is to iterate the previous procedure with $g_1$ instead of $g$ such that the unwanted term in the decomposition (\ref{eqn:g1}) becomes smaller and smaller at each step. We shall see that if $m$ is chosen big enough and fixed, this error term disappears as the number of iterations tends to $+\infty$.

\noindent If we replace $g$ by $g_1$, we get the following decomposition analogous to (\ref{eqn:g1}):  \\

$\displaystyle g_1(z)=\sum\limits_{j=1}^{N_m} a_{1, \, j} \, \sigma_{m, \, j}(z) + \sum\limits_{l=1}^{n}z_l \, v_{1, \, l}(z) + g_2(z) ,$ \hspace{2ex} $z\in B,$ \hspace{2ex} with \\

  $\displaystyle \int_B |g_2|^2 \, e^{-2m\varphi} \leq \bigg(C(r, \, r', \, d)\, \varepsilon\bigg)^2 \, C_g,$ \hspace{3ex}  \\

 $\displaystyle \sum\limits_{j=1}^{N_m}|a_{1, \, j}|^2 \leq  C(r, \, r', \, d) \, \bigg(C(r, \, r', \, d)\, \varepsilon\bigg) \, C_g,$ \\

  $\displaystyle \sum_{l=1}^{n} \int_B |v_{1, \,l}|^2 \, e^{-2m\varphi} \leq C_1(r, \, r', \, d)\, \bigg(C(r, \, r', \, d)\, \varepsilon\bigg)\, C_g,$  \\

\noindent which further gives:  \\

$\displaystyle g(z) = \sum\limits_{j=1}^{N_m} \bigg(a_j + a_{1, \, j}\bigg)\, \sigma_{m, \, j}(z) + \sum\limits_{l=1}^{n}z_l \bigg(v_l(z) + v_{1, \, l}(z)\bigg)+ g_2(z),$  \\

 \noindent for every $z\in B.$ We can now indefinitely iterate this procedure with $g_2$ instead of the original $g$. After $p$ iterations, we get the decomposition :  \\

$\displaystyle g(z) = \sum\limits_{j=1}^{N_m} \bigg(a_j + a_{1, \, j} + \dots + a_{p-1, \, j}\bigg)\, \sigma_{m, \, j}(z) +$  \\

\hspace{15ex} $ + \sum\limits_{l=1}^{n}z_l \bigg(v_l(z) + v_{1, \, l}(z) + \dots + v_{p-1, \, l}(z)\bigg)+ g_p(z),$  \\

\noindent for every $z\in B,$ with estimates:

\begin{equation}\label{eqn:gps}\displaystyle \int_B |g_p|^2 \, e^{-2m\varphi} \leq \bigg(C(r, \, r', \, d)\, \varepsilon\bigg)^p \, C_g,\end{equation}   

 $\displaystyle \sum\limits_{j=1}^{N_m}|a_{p-1, \, j}|^2 \leq  C(r, \, r', \, d) \, \bigg(C(r, \, r', \, d)\, \varepsilon\bigg)^{p-1} \, C_g,$ \\

  $\displaystyle \sum_{l=1}^{n} \int_B |v_{p-1, \,l}|^2 \, e^{-2m\varphi} \leq C_1(r, \, r', \, d)\, \bigg(C(r, \, r', \, d)\, \varepsilon\bigg)^{p-1}\, C_g,$  \\

\noindent Let us now set $\varepsilon_0:= C(r, \, r', \, d)\, \varepsilon,$ $c_j^{(p)}:=a_j + a_{1, \, j} + \dots + a_{p-1, \, j},$ and $v_{l}^{(p)}:=v_l + v_{1, \, l} + \dots + v_{p-1, \, l}.$ On the one hand, we get:  \\

$\displaystyle \bigg(\sum\limits_{j=1}^{N_m}|c_j^{(p)}|^2\bigg)^{\frac{1}{2}} \leq \sum\limits_{j=1}^{N_m}|c_j^{(p)}| \leq \sum\limits_{j=1}^{N_m}|a_j| + \sum\limits_{j=1}^{N_m}|a_{1, \, j}| + \dots + \sum\limits_{j=1}^{N_m} |a_{p-1, \, j}| $  \\

\noindent $ \leq \bigg(N_m \, C(r, \, r', \, d)\, C_g\bigg)^{\frac{1}{2}}\, (1 + \varepsilon_0^{\frac{1}{2}} + \dots + \varepsilon_0^{\frac{p-1}{2}}) = \bigg(N_m \, C(r, \, r', \, d)\, C_g\bigg)^{\frac{1}{2}}\, \frac{1-\varepsilon_0^{\frac{p}{2}}}{1-\varepsilon_0^{\frac{1}{2}}}.$  \\

\noindent On the other hand, if $||\cdot ||$ stands for the $e^{-2m\varphi}$-weighted $L^2$ norm, we get: \\

$\sum\limits_{l=1}^{n}||v_l^{(p)}||\leq \sqrt{n} \, \bigg[ (\sum\limits_{l=1}^{n}||v_l||^2)^{\frac{1}{2}} +(\sum\limits_{l=1}^{n}||v_{1, \,l}||^2)^{\frac{1}{2}} + \dots + (\sum\limits_{l=1}^{n}||v_{p-1, \,l}||^2)^{\frac{1}{2}}\bigg]$  \\

$\displaystyle \leq \sqrt{n\, C_1(r, \, r', \, d)\, C_g}\, (1+ \varepsilon_0^{\frac{1}{2}}+ \dots + \varepsilon_0^{\frac{p-1}{2}})= \sqrt{n\, C_1(r, \, r', \, d)\, C_g}\, \frac{1-\varepsilon_0^{\frac{p}{2}}}{1-\sqrt{\varepsilon_0}},$ \\

\noindent by using the Cauchy-Schwarz inequality and the previous estimates on $v_{k, \, l}$'s. We finally get:  \\

$\displaystyle\sum\limits_{l=1}^{n}||v_l^{(p)}||^2 \leq \bigg(\sum\limits_{l=1}^{n}||v_l^{(p)}||\bigg)^2 \leq n\, C_1(r, \, r', \, d)\, \bigg(\frac{1-\varepsilon_0^{\frac{p}{2}}}{1-\sqrt{\varepsilon_0}}\bigg)^2 \, C_g.$ \\

\noindent for every $p\in \mathbb{N}^{\star}.$ Let us now choose $\varepsilon$ so small that $\varepsilon_0:= C(r, \, r', \, d)\, \varepsilon =\frac{1}{4}.$ By letting $p\rightarrow +\infty$, the estimate (\ref{eqn:gps}) shows that $g_p$ disappears from the expression of $g$, and we get:  \\

$\displaystyle g(z) = \sum\limits_{j=1}^{N_m} c_j \, \sigma_{m, \, j}(z) + \sum\limits_{l=1}^{n}z_l h_l(z)$,  \hspace{3ex} $z\in B$,  \\

\noindent with $c_j\in \C$ obtained as limits of $c_j^{(p)}$, and holomorphic functions $h_l$ on $B$ obtained as limits of $v_l^{(p)}$, satisfying the estimates in the statement, after possibly replacing $4\, C(r, \, r', d)$ and $4\, n\, C_1(r, \, r', \, d)$ by the maximum $C$ of these two constants. Lemma \ref{Lem:firstordergeneration} is proved. \hfill $\Box$

\vspace{2ex}

 We can now run an induction argument using lemma \ref{Lem:firstordergeneration} repeatedly to get, at every step $p$, an approximation to order $p$ of the original local section $g$ of ${\cal I}(m\varphi)$ by a finite linear combination of $\sigma_{m, \, j}$'s. The following is a slightly more precise rewording of the first part of Theorem \ref{The:introd1}.

\begin{The}\label{The:finitegeneration} Given $\varphi$ such that $i\partial\bar{\partial}\varphi \geq C_0 \, \omega$, and having fixed a ball $B:=B(x, \, r)\subset\subset \Omega$, there exist a ball $B(x, \, r_0)\subset\subset B(x, \, r),$ $\varepsilon>0,$ and $m_0=m_0(C_0)\in \mathbb{N}$ such that for every $m\geq m_0$ the following property holds. Every $g\in {\cal H}_B(m\varphi)$ admits, with respect to the $\varepsilon$-concentrated $\sigma_{m, \, j}$'s in an orthonormal basis of ${\cal H}_{\Omega}(m\varphi)$ made up of eigenvectors of $T_{B, \, m}$, a decomposition:  \\

$\displaystyle g(z) = \sum\limits_{j=1}^{N_m} b_{m, \, j}(z)\, \sigma_{m, \, j}(z)$,  \hspace{3ex}  $z\in B(x, \, r_0),$  \\

\noindent with some holomorphic functions $b_{m, \, j}$ on $B(x, \, r_0),$ satisfying: \\

 $\displaystyle\sup\limits_{B(x, \, r_0)}\sum\limits_{j=1}^{N_m}|b_{m, \, j}|^2 \leq C \, N_m \, \int_B|g|^2\, e^{-2m\varphi}<+\infty,$ \\

\noindent where $C>0$ is a constant depending only on $n, \, r,$ and the diameter of $\Omega.$

\end{The} 

\vspace{1ex}

\noindent {\it Proof.} We may assume, without loss generality, that $x=0$. Fix $\varepsilon>0$ as in lemma \ref{Lem:firstordergeneration}, associated with an arbitrary fixed open ball $B(0, \, r') \subset\subset B(0, \, r)$. Let $g\in {\cal H}_B(m\varphi)$ with $C_g:=\int_B|g|^2\, e^{-2m\varphi}<+\infty$. Lemma \ref{Lem:firstordergeneration} gives a decomposition $g(z) = \sum\limits_{j=1}^{N_m} a_j\, \sigma_{m, \, j}(z) + \sum\limits_{l_1=1}^{n} z_{l_1} \, h_{l_1}(z)$ for every $z\in B(0, \, r)$, with coefficients under control. We now apply lemma \ref{Lem:firstordergeneration} to every function $h_{l_1}$, $l_1=1, \dots , n,$ and get:  \\

$\displaystyle h_{l_1}(z)= \sum\limits_{j=1}^{N_m} a_{j, \, l_1}\, \sigma_{m, \, j}(z) + \sum\limits_{l_2=1}^{n} z_{l_2} \, h_{l_1, \, l_2}(z)$, \hspace{3ex} $z\in B(0, \, r)$, \\

\noindent with constant coefficients $a_{j, \, l_1} \in \C$, satisfying: \\

$\displaystyle \sum\limits_{l_1=1}^{n} \sum\limits_{j=1}^{N_m}|a_{j, \, l_1}|^2 \leq C\, N_m \, \sum\limits_{l_1=1}^{n} \int_B|h_{l_1}|^2 \, e^{-2m\varphi} \leq C^2 \, N_m \, C_g,$  \\

\noindent and holomorphic functions $h_{l_1, \, l_2}$ on $B$, satisfying:  \\

$\displaystyle \sum\limits_{l_1=1}^{n} \sum\limits_{l_2=1}^{n} \int_B|h_{l_1, \, l_2}|^2 \, e^{-2m\varphi} \leq C \, \sum\limits_{l_1=1}^{n} \int_B|h_{l_1}|^2 \, e^{-2m\varphi} \leq C^2 \, C_g.$  \\

\noindent We thus obtain, after $p$ applications of lemma \ref{Lem:firstordergeneration}, the decomposition:  \\

\noindent $\displaystyle g(z)=\sum\limits_{j=1}^{N_m} \bigg(a_j + \sum\limits_{\nu=1}^{p-1}\sum\limits_{l_1, \dots , \, l_{\nu}=1}^{n} a_{j, \, l_1, \dots , \, l_{\nu}}\, z_{l_1} \dots z_{l_{\nu}}\bigg) \, \sigma_{m, \, j}(z) + \sum\limits_{l_1, \dots , \, l_p=1}^{n} \, z_{l_1} \dots z_{l_p}\, v_{l_1, \dots , \, l_p}(z),$  \\

\noindent for every $z\in B(0, \, r)$, with coefficients $a_{j, \, l_1, \dots , \, l_{\nu}=1}\in \C$ and $v_{l_1, \dots , \, l_p}\in {\cal O}(B(0,\, r))$, satisfying the estimates:  \\

$\displaystyle \sum\limits_{l_1, \dots , \, l_{\nu}=1}^{n}\sum\limits_{j=1}^{N_m}|a_{j, \, l_1, \dots , \, l_{\nu}}|^2 \leq  C^{\nu +1}\, N_m \, C_g,$ \hspace{3ex} $\nu = 1, \dots , \, p-1,$  \\

\noindent $\displaystyle \sum\limits_{l_1, \dots , \, l_p=1}^{n} \int_B |v_{l_1, \dots , \, l_p}|^2 \, e^{-2m\varphi} \leq  C^p \, C_g.$  \\

\noindent We now put $\displaystyle b_{m, \, j}(z):=a_j + \sum\limits_{\nu=1}^{+\infty}\sum\limits_{l_1, \dots , \, l_{\nu}=1}^{n} a_{j, \, l_1, \dots , \, l_{\nu}}\, z_{l_1} \dots z_{l_{\nu}},$ for $j=1, \dots , N_m$, and will prove that the series defining $b_{m, \, j}$ converges to a holomorphic function on some smaller ball $B_0:=B(0, \, r_0),$ and that, with an appropriate choice of $r_0$, we have:  \\

$\displaystyle \sup_{B(0, \, r_0)}\sum_{j=1}^{N_m} |b_{m, \, j}|^2 \leq \frac{1}{(1-r/d)^2} \, C\, N_m \, C_g.$  \\

\noindent Since $\sup\limits_{B(0, \, r_0)} |z_{l_1} \dots z_{l_{\nu}}|^2 \leq r_0^{2\nu},$ we get, for every $1 \leq \nu < +\infty$ and every $z\in B(0, \, r_0),$ the estimate:  \\

\vspace{2ex}

$\displaystyle \sum\limits_{j=1}^{N_m} \bigg|\sum\limits_{l_1, \dots , \, l_{\nu}=1}^{n} a_{j, \, l_1, \dots , \, l_{\nu}}\, z_{l_1} \dots z_{l_{\nu}}\bigg|^2 \leq r_0^{2\nu} \, \sum\limits_{j=1}^{N_m} \bigg(\sum\limits_{l_1, \dots , \, l_{\nu}=1}^{n} |a_{j, \, l_1, \dots , \, l_{\nu}}|\bigg)^2$  \\

\hspace{3ex} $\displaystyle\leq r_0^{2\nu} \, n^{\nu} \, \sum\limits_{j=1}^{N_m}\sum\limits_{l_1, \dots , \, l_{\nu}=1}^{n}|a_{j, \, l_1, \dots , \, l_{\nu}}|^2 \leq \bigg(n\, r_0^2\, C \bigg)^{\nu} \, C\, N_m \, C_g$ \\

\hspace{3ex} $\displaystyle = \bigg(\frac{r}{d}\bigg)^{2\nu}\, C\, N_m\, C_g,$  \\

\vspace{2ex}

\noindent if we choose $\displaystyle r_0=\frac{r}{d} \, \frac{1}{\sqrt{n\, C}}.$ The remaining arguments are purely formal. Put:  \\

$F_{\nu, \, j}:=\sum\limits_{l_1, \dots , \, l_{\nu}=1}^{n} a_{j, \, l_1, \dots , \, l_{\nu}}\, z_{l_1} \dots z_{l_{\nu}},$ \hspace{2ex} for $\nu\geq 1$ and $j=1, \dots , \, N_m,$  \\

$F_{0, \, j}=a_j$, \hspace{2ex} for $j=1, \dots , \, N_m.$  \\

\noindent Since $\displaystyle b_{m, \, j}=\sum\limits_{\nu=0}^{+\infty} F_{\nu, \, j},$ for $j=1, \dots , \, N_m,$ we get:  \\

$\displaystyle \bigg(\sum_{j=1}^{N_m}|b_{m, \, j}|^2\bigg)^{\frac{1}{2}} \leq \sum\limits_{\nu=0}^{+\infty}\bigg(\sum_{j=1}^{N_m}|F_{\nu, \, j}|^2\bigg)^{\frac{1}{2}} \leq \sum\limits_{\nu=0}^{+\infty}\bigg(\frac{r}{d}\bigg)^{\nu}\, \sqrt{C\, N_m\, C_g}$  \\

\hfill $\displaystyle = \frac{1}{1-r/d}\, \sqrt{C\, N_m\, C_g},$  \\

\noindent at every point in $B(0, \, r_0).$ If we absorb the denominator in the constant $C>0$, the proof is complete. \hfill $\Box$

\begin{Rem}\label{Rem:onintrod} In retrospect, we see that estimate (\ref{eqn:numberestimate}) which concluded section \ref{subsection:Toeplitz} clearly proves the last statement of Theorem \ref{The:introd1} in the introduction. The proof of this theorem is now complete.

\end{Rem}

\begin{Cor}\label{Cor:finitegeneration} Under the hypotheses of Theorem \ref{The:finitegeneration}, the following estimate holds:  \\

$\displaystyle \sum\limits_{j=1}^{+\infty}|\sigma_{m, \, j}(z)|^2 \leq C\, N_m \,  \sum\limits_{j=1}^{N_m}|\sigma_{m, \, j}(z)|^2$, \hspace{3ex} $z\in B(x, \, r_0)$.

\end{Cor}

\noindent {\it Proof.} As $\displaystyle \sum\limits_{j=1}^{+\infty}|\sigma_{m, \, j}(z)|^2 = \sup\limits_{f\in \bar{B}_m(1)} |f(z)|^2$, where $\bar{B}_m(1)$ is the closed unit ball in ${\cal H}_{\Omega}(m\varphi)$, and as every such $f$ has a decomposition as in Theorem \ref{The:finitegeneration} on $B(x, \, r_0)$, the estimate follows from the Cauchy-Schwarz inequality. \hfill $\Box$

\vspace{3ex}

\subsection {Regularization of plurisubharmonic functions}\label{subsection:regularization}

  We now turn to the second part of the paper. A by now classical result of Demailly's ([Dem92], {\it Proposition 3.1.}) states that a psh function $\varphi$ with arbitrary singularities can be approximated by psh functions $\varphi_m$ with analytic singularities (cf. $(\star)$), constructed as  

 \begin{equation}\label{eqn:defreg}\displaystyle\varphi_m(z) = \frac{1}{2m} \, \log\sum\limits_{j=1}^{+\infty}|\sigma_{m, \, j}(z)|^2= \sup\limits_{f\in \bar{B}_m(1)}\, \frac{1}{m}\, \log |f(z)|,  \hspace{3ex} z\in \Omega,\end{equation}

\noindent where $(\sigma_{m, \, j})_{j\in \mathbb{N}^{\star}}$ is an arbitrary orthonormal basis, and $\bar{B}_m(1)$ is the unit ball, of the Hilbert space ${\cal H}_{\Omega}(m\varphi)$ considered in the previous sections. We even have:

  \begin{equation}\label{eqn:approximation}\varphi(z) - \frac{C_1}{m} \leq \varphi_m(z) \leq \sup_{|\zeta - z|<r}\varphi(\zeta) + \frac{1}{m} \, \log(\frac{C_2}{r^n}),\end{equation}

\noindent for every $z\in \Omega$ and every $r< d(z, \, \partial \Omega)$. The lower estimate is a consequence of the Ohsawa-Takegoshi $L^2$ extension theorem. The upper estimate is far easier, coming from an application of the submean value inequality satisfied by squares of absolute values of holomorphic functions. Our aim is to improve the upper estimate by replacing the supremum by a pointwise upper bound which affords a much better understanding of singularities. This is not possible for an arbitrary $\varphi$, but the following proposition shows it to be possible if $\varphi$ is assumed to have analytic singularities.

 \begin{The}\label{The:upperbound} Let $\Omega\subset\subset \mathbb{C}^n$ be a bounded pseudoconvex open set, and let $\varphi=\frac{c}{2} \, \log (|g_1|^2 + \dots + |g_N|^2)$ be a plurisubharmonic function with analytic singularities on $\Omega$. If $\Omega' \subset\subset \Omega'' \subset\subset \Omega$ are relatively compact open subsets, then for every $\delta >0$ and every $m\geq \frac{n+2}{c\, \delta}$, we have: \\

$\displaystyle \varphi_m(z) \leq (1-\delta)\, \varphi(z) + c\, \delta \, \log A + \frac{\log(C_n\, (m\, c -n))}{m}$, \hspace{2ex} $z\in \Omega'$, \\

\noindent where $A:=\max\{\sup\limits_{\Omega''}\sqrt{|g_1|^2 + \dots + |g_N|^2}, \, 1\},$ and $C_n>0$ is a constant depending only on $\Omega', \, \Omega''$, and $n$.

\end{The}

\begin{Rem} {\rm The upper bound given in Theorem \ref{The:upperbound}, combined with Demailly's lower bound in inequality \ref{eqn:approximation}, implies that $\varphi$ and its regularizations $\varphi_m$ have the same $-\infty$ poles on $\Omega'$, if $m$ is big enough. Only the Lelong numbers may be slightly different. Since the set of poles of an arbitrary psh function  $\varphi$ is not necessarily analytic, whereas the polar set of $\varphi_m$ is always analytic, we see that some analyticity assumption on $\varphi$ is necessary. We may ask for the weakest such assumption on $\varphi$ under which Theorem \ref{The:upperbound} holds.}

\end{Rem}

\noindent {\it Proof.} Fix $0<\delta <1$. Let $m\geq \frac{n+2}{c\, \delta}$, and let $f$ be an arbitrary holomorphic function on $\Omega$ such that: \\

 $\displaystyle \int_{\Omega}|f|^2 \, e^{-2m\varphi}=1 \hspace{2ex} \Leftrightarrow \hspace{2ex} \int_{\Omega}\frac{|f|^2}{(|g_1|^2 + \dots + |g_N|^2)^{mc}} =1$. \\

\noindent  We clearly have $mc \geq n+2$; put $q:=[mc]-(n+1)\geq 1$. Then Skoda's $L^2$ division theorem ([Sko72]) gives the existence of holomorphic functions $h_{i_1, \dots , i_q}$ on $\Omega$, for all multiindices $(i_1, \dots, i_q)\in \{1, \dots , N\}^q$, such that:  \\

$\displaystyle f(z)= \sum\limits_{i_1, \dots , i_q=1}^{N} h_{i_1, \dots , i_q}(z) \, g_{i_1}(z) \dots g_{i_q}(z)$, \hspace{2ex} $z\in \Omega$,  \\

\noindent and \\

$\displaystyle \sum\limits_{i_1, \dots , i_q=1}^{N} \int_{\Omega} \frac{|h_{i_1, \dots , i_q}|^2}{(|g_1|^2 + \dots + |g_N|^2)^{mc-q}} \leq $ \\

\hspace{10ex} $\displaystyle \leq \frac{m \, c-n}{m\, c-[m\, c]+1}\, \int_{\Omega}\frac{|f|^2}{(|g_1|^2 + \dots + |g_N|^2)^{mc}} \leq m\, c-n$.  \\    

\noindent Put $|g|^2:= |g_1|^2 + \dots + |g_N|^2.$ As $\Omega''$ satisfies $\Omega' \subset\subset \Omega'' \subset\subset \Omega$, we get, for every $z\in \Omega'$:  \\

$\displaystyle  |f(z)|^2 \leq |g(z)|^{2q} \, \sum\limits_{i_1, \dots , i_q=1}^{N}|h_{i_1, \dots , i_q}(z)|^2 \leq    $  \\

\hspace{6ex} $\displaystyle \leq C_n \, |g(z)|^{2q} \, \sum\limits_{i_1, \dots , i_q=1}^{N} \int\limits_{\Omega''} \frac{|h_{i_1, \dots , i_q}|^2}{(|g_1|^2 + \dots + |g_N|^2)^{mc-q}}\, |g|^{2(mc-q)}$  \\

\hspace{6ex} $\displaystyle \leq C_n \, (m\, c-n) \, (\sup\limits_{\Omega''}|g|)^{2(mc-q)}\, |g(z)|^{2q},$  \\

\noindent where $C_n>0$ is the constant involved in the submean value inequality applied to every $|h_{i_1, \dots , i_q}|^2$ at $z$, and thus depending only on the open sets $\Omega'$, $\Omega''$, and on $n$, but independent of $m$. We can see, since $A=\max\{\sup\limits_{\Omega''}|g|, \, 1\},$ that:  \\

$C_n \, (m\, c-n) \, (\sup\limits_{\Omega''}|g|)^{2(mc-q)} \leq  C_n\, (m\, c-n)  \, A^{2(mc-q)}.$ \\

\noindent  We thus get: \\

$\displaystyle |f(z)|^2 \, e^{-2m(1-\delta)\varphi(z)} \leq C_n\, (m\, c-n)\, A^{2(mc-q)}   \, \frac{|g(z)|^{2([mc]-n-1)}}{|g(z)|^{2mc(1-\delta)}}$ \\

\hspace{10ex} $\displaystyle = C_n\, (m\, c-n)\, A^{2(mc-q)} \, |g(z)|^{2(mc\delta - n -1-(mc - [mc]))},$ \hspace{3ex} $z\in \Omega'$.  \\

\noindent Since $m$ has been chosen such that $mc\delta - n- 2\geq 0$, the exponent $2(mc\delta - n -1-(mc - [mc]))$ above is positive, and we thus finally get: \\

$\displaystyle |f(z)|^2 \, e^{-2m(1-\delta)\varphi(z)} \leq C_n\, (m\, c-n)\, A^{2(mc-q)} \, A^{2(mc\delta - n -1-(mc - [mc]))}$  \\

 \hspace{8ex} $= C_n\, (m\, c-n)\, A^{2(mc\delta - n -1 +[mc]-q)} =C_n\, (m\, c-n)\, A^{2mc\delta},$  \hspace{2ex} $z\in \Omega'$.  \\

\noindent By taking $\log$  and dividing by $2m$ on both sides, we get: \\

$\displaystyle \varphi_m(z) \leq (1-\delta)\, \varphi(z) + c\delta \, \log A + \frac{\log(C_n(m\, c-n))}{m}$, \hspace{2ex} $z\in \Omega'$,  \\

\noindent  The proof is complete.  \hfill $\Box$

\begin{Rem}\label{Rem:upperbound} If in Theorem \ref{The:upperbound} the psh function with analytic singularities is given, more generally, as $\varphi = \frac{c}{2}\, \log(|g_1|^2 + \dots + |g_N|^2) + v,$ for some $C^{\infty}$ psh function $v,$ then the estimate is modified as:  \\

$\displaystyle \varphi_m(z) \leq (1-\delta)\, \varphi(z) + c\,\delta \, \log A + \frac{\log(C(m\, c-n))}{m} + \sup\limits_{\Omega''}v-(1-\delta)\, \inf\limits_{\Omega''}v$,   \\

\noindent for every $z\in \Omega'$ and every $m\geq \frac{n+2}{c\, \delta}.$

\end{Rem}

\vspace{3ex}

Theorem \ref{The:introd2} stated in the Introduction appears now as a consequence of Theorem \ref{The:upperbound}.  \\

\noindent {\bf Proof of Theorem \ref{The:introd2}.} As repeatedly pointed out above, the ideal sheaf ${\cal I}(m\varphi)$ is generated as an ${\cal O}_{\Omega}$-module by the Hilbert space ${\cal H}_{\Omega}(m\varphi)$. Now, by the definition (\ref{eqn:defreg}) of $\varphi_m$, the estimate obtained in Theorem \ref{The:upperbound} implies that every $f\in {\cal H}_{\Omega}(m\varphi)$ with norm $1$ satisfies: \\

\hspace{6ex}$\displaystyle |f(z)|^2  \leq C\, (m\, c-n)\, A^{2\,m\,c\,\delta}\, e^{2m(1-\delta)\varphi(z)},$ \hspace{3ex} $z\in \Omega',$ $m\geq \frac{n+2}{c\, \delta}.$ \\

\noindent  If $f_1, \dots , \, f_p\in {\cal H}_{\Omega}(m\varphi)$ have norm $1,$ then $|f_1 \dots f_p|^2 \, e^{-2mp(1-\delta)\varphi}$ is bounded and therefore integrable on $\Omega'$, proving that $f_1, \dots , \, f_p$ is a section of ${\cal I}(mp(1-\delta)\varphi)$ on $\Omega'.$ The proof is complete.  \hfill $\Box$

\vspace{5ex}

\noindent {\bf\Large References}  \\

\noindent [Ber03] \, B. Berndtsson --- {\it Bergman Kernels Related to Hermitian Line Bundles Over Compact Complex Manifolds} ---  Explorations in complex and Riemannian geometry,  1--17, Contemp. Math., 332, Amer. Math. Soc., Providence, RI, 2003.

\vspace{1ex}

\noindent [Dem92] \, J.-P. Demailly --- {\it Regularization of Closed Positive Currents and Intersection Theory} --- J. Alg. Geom., {\bf 1} (1992), 361-409.

\vspace{1ex}

\noindent [DEL00] \, J.-P. Demailly, L. Ein, R. Lazarsfeld --- {\it A Subadditivity Property of Multiplier Ideals}--- Michigan Math. J., {\bf 48} (2000), 137-156.

\vspace{1ex}

\noindent [H\"or65] \, L. H\"ormander --- {\it $L^2$ Estimates and Existence Theorems for the $\bar{\partial}$ Operator} --- Acta Math. {\bf 113} (1965), 89-152.

\vspace{1ex}

\noindent [Lan67] \, H. J. Landau --- {\it Necessary Density Conditions for Sampling and Interpolation of Certain Entire Functions} --- Acta Math. {\bf 117} (1967), 37-52.

\vspace{1ex}

\noindent [Lin01] \, N. Lindholm --- {\it Sampling in Weighted $L\sp p$ Spaces of Entire Functions in ${\Bbb C}\sp n$ and Estimates of the Bergman Kernel} ---  J. Funct. Anal.  18  (2001),  no. 2, 390--426.

\vspace{1ex}

\noindent [Siu02] \, Y-T.Siu --- {\it Extension of Twisted Pluricanonical Sections with Plurisubharmonic Weight and Invariance of Semipositively Twisted Plurigenera for Manifolds Not Necessarily of General Type} --- Complex geometry (G\"ottingen, 2000),  223--277, Springer, Berlin, 2002.

\vspace{1ex}

\noindent [Sko72] \, H. Skoda --- {\it Applications des techniques $L\sp{2}$ \`a la th\'eorie des id\'eaux d'une alg\`ebre de fonctions holomorphes avec poids}---  Ann. Sci. \'Ecole Norm. Sup. (4)  5  (1972), 545--579.

\vspace{6ex}

\noindent Mathematics Institute, University of Warwick, Coventry CV4 7AL, UK 

\noindent E-mail: popovici@maths.warwick.ac.uk

\end{document}